\documentclass[a4paper]{article}

\usepackage[utf8x]{inputenc}
\usepackage{amssymb,amsmath,amsthm}
\usepackage{graphicx}
\usepackage{subfigure}
\usepackage{color}
\usepackage{psfrag}
\usepackage{relsize}

\newtheorem{mydef}{Definition}
\newtheorem{mytheorem}{Theorem}

\newtheorem{myalgorithm}{Algorithm}

\DeclareMathOperator*{\Prob}{Prob}

\newtheorem*{mySTOP}{Algorithmic Stop Criterion}

\newcommand{\corrFOURDOM}{\textcolor{black}}{}
\newcommand{\corrFIVEDOM}{\textcolor{black}}{}

\begin{document}

\title{A robust optimization approach to experimental design for model
  discrimination of dynamical systems}
\author{Dominik Skanda$^1$ \and Dirk Lebiedz$^2$}

\maketitle

\footnotetext[1]{  University of Freiburg \\
	   Center for Systems Biology (ZBSA) \\
           Habsburgerstr. 49 \\
	   79104 Freiburg \\
	   Germany\\
           $^2$ University of Freiburg \\
	   Center for Systems Biology (ZBSA) \\
	   Faculty of Biology \\
	   and Department of Applied Mathematics, \\
	   Habsburgerstr. 49 \\
	   79104 Freiburg \\
	   Germany \\
	} 

\maketitle

\begin{abstract}
A high-ranking goal of interdisciplinary modeling approaches in science and
engineering are quantitative prediction of system dynamics and 
model based optimization. Quantitative modeling has to be closely related to experimental investigations if the model is supposed to be used for mechanistic analysis and model predictions. Typically, before an appropriate model of an experimental system is found different hypothetical models might be reasonable and consistent with previous knowledge and available data. 
The parameters of the models up to an estimated confidence region are
generally not known a priori. Therefore one has to incorporate possible
parameter configurations of different models into a model discrimination
algorithm which leads to the need for robustification. 
In this article we present a numerical algorithm which calculates a design of
experiments allowing optimal discrimination of different hypothetic candidate
models of a given dynamical system for the most inappropriate (worst case)
parameter configurations within a parameter range. The design
comprises initial values, system perturbations and the optimal
\corrFIVEDOM{placement of} measurement time points, the number of measurements
as well as the time \corrFIVEDOM{points} are subject to 
design. The statistical discrimination criterion is worked out rigorously for
these settings, a derivation from the Kullback-Leibler divergence as
optimization objective is presented for the case of discontinuous
Heaviside-functions modeling the measurement decision which are replaced by
continuous approximation during the 
optimization procedure. The resulting problem can be classified as a semi-infinite optimization problem which we solve in an outer approximations
approach stabilized by a suggested homotopy strategy whose efficiency is
demonstrated. We present the theoretical framework, algorithmic
realization \corrFIVEDOM{ and numerical results}. 
\end{abstract}

\section{Introduction}

High-ranking goals of interdisciplinary modeling approaches in the natural
sciences are quantitative prediction of system dynamics and 
model based optimization. In particular in modern systems biology a
related issue is to link molecular attributes to dynamic mechanisms and functional
properties at the system level in order to mechanistically understand emerging
functionality. For these purposes, mathematical modeling,
numerical simulation and scientific computing techniques are indispensable.
Quantitative modeling closely combined with experimental investigations is
required if the model is supposed to be used for sound mechanistic
analysis and model predictions.  

Typically, before an appropriate model of a system is found different
hypothetical models might be reasonable and consistent with previous knowledge
and available data. 
The goal of this article is to derive, develop, implement and apply a
numerical algorithm which calculates in a suitable sense an optimal design of
experiments which allows the best 
discrimination of different hypothetic candidate models in \corrFOURDOM{form} of ordinary
differential equations (ODE). The algorithmic idea is to iteratively separate
the response of different models by use of variations of
experimental conditions and perturbations to the system.  

To discriminate a set of candidate models against a given set of experimental
data likelihood ratio tests based on bootstrap methods have been described in
the literature, see e.g. \cite{Horn1987}, \cite{Stricker1994} or \cite{Timmer2004}. 
Ranking methods like Stewart's method (\cite{Stewart1998}) or the well known
Akaike information criterion (see e.g. \cite{Burnham2002}) are popular as well
in the field of biological modeling.  
Applications can be found for example in \cite{RishiJain2006} or \cite{Bernacki2009}.\\ 
In contrast to these approaches our work deals
with the problem of designing experiments so that statistical methods can be
exploited in an optimal sense for model discrimination. This differs from the approach to find an
experimental design to best estimate the parameters of a model for a given
experimental system in terms of criteria characterizing the confidence regions
\cite{Korkel1999,Bauer2000,Balsao2008}. 

Different approaches to design experiments for model discrimination
exist. Besides optimization methods (see e.g \cite{Lacey1984},
\cite{Cooney1995}, \cite{Takors1997} or \cite{Kremling2004}) a model-based
feedback controller see e.g \cite{Apgar2008} and Markov chain Monte Carlo
sampling methods \cite{Myung2009} have been used to construct an appropriate
design. An overview of \corrFIVEDOM{various} experimental design techniques can be found
in \cite{Kreutz2009}.  \\
 \corrFIVEDOM{Here,} we present a robust numerical optimization algorithm which calculates the
 optimal design of experiments allowing the best discrimination of different
 candidate ODE models. An appropriate model and its parameters up to an estimated
 confidence region are not known a priori. Therefore one has to incorporate
 possible parameter configurations of different models into a model
 discrimination algorithm. The aim is to calculate the most discriminable
 response of different models for the most inappropriate parameter
 configurations within a parameter range via a worst case estimate. In that
 context inappropriate parameter configurations refer to the case when
 different models have calibrated parameter values such that the models have
 the most similar response. This worst case estimate leads to the formulation
 of a maxmin optimization problem. 
Building on our previous work \cite{Skanda2010} we present an algorithm to
compute robust optimal experimental designs. For the robustification we set up
an outer approximation approach stabilized by a homotopy strategy.\\ 
The article is organized as follows. In Section \ref{kloptimal} we give a
brief overview of so called Kullback-Leibler(KL)-optimality as discussed by
L\'opez-Fidalgo et al. \cite{Lopez-Fidalgo2007} in the context of model
discrimination. In Section \ref{biodesign} we derive our optimal design criterion
by use of KL-optimality. In Section 
\ref{maxminTheory} the theoretical framework for the calculation of a robust
design via solution of a maxmin optimization problem is presented. The
numerical implementation is discussed in Section \ref{numericRealization}. The
homotopy solution strategy is presented in Section \ref{HomotopyStrategy}. We
demonstrate applications of the algorithmic framework for two test cases from
biology, an allosteric metabolic enzyme model for glycolytic oscillations and
a model describing signal sensing in dictyostelium discoideum, results are
presented in Section \ref{numericalResults}. 
\section{Statistical Basis of Model Discrimination}
\label{theory}
\subsection{KL-optimal design}
\label{kloptimal}
In this section a model discrimination criterion based on the Kullback-Leibler
(KL) divergence called KL-optimality as discussed by L\'opez-Fidalgo et al. \cite{Lopez-Fidalgo2007} is
introduced. L\'opez-Fidalgo et al.  \cite{Lopez-Fidalgo2007} demonstrate that KL-optimality
is consistent with T-optimality \cite{ATKINSON1975} and generalized
T-optimality \cite{Ucinski2004} which are well known model discrimination
criteria based on statistical tests.\\\\
We introduce the concept of a probability space and formally define the KL-divergence. 
\begin{mydef}
 The $probability$ $space$ is a triple $\begin{pmatrix}\Omega,\mathcal{F},P \end{pmatrix}$ consisting of
 \begin{itemize}
  \item a non-empty set $\Omega$ (sample space),
  \item a $\sigma$-algebra $\mathcal{F}\subseteq \mathcal{P}(\Omega)$, $E \in
    \mathcal{F}$ is called an event
  \item a probability measure $P:\mathcal{F}\rightarrow[0,1]$.
 \end{itemize}
\end{mydef}

\begin{mydef}
 Two probability spaces $\begin{pmatrix}\Omega,\mathcal{F},P_i \end{pmatrix}$,
 $i=1,2$, are called $absolutely$ $continuous$ with respect to each other, in
 symbols $P_1 \equiv P_2$, if $\nexists \ E\in\mathcal{F}: (P_1(E)=0$ AND $P_2(E)\neq0)$ OR $(P_1(E)\neq0$ AND $P_2(E)=0)$. 
\end{mydef}

The Radon-Nikodym Theorem allows a representation of a probability measure via
a measurable probability density function.

\begin{mytheorem}\label{theorem:radon_nikodym}
(Radon-Nikodym)\\
Let $\lambda$ be a probability measure such that $\lambda\equiv P_1$,
$\lambda\equiv P_2$. Then $\lambda$-measurable functions $f_i:\Omega
\rightarrow \mathbb{R}, i=1,2$,
called $generalized$ $probability$ $densities$, exist which are unique up to
sets of measure zero and non-negative, such that 
\begin{equation}
 P_i(E)=\int_Ef_i(x)\mathrm{d}\lambda(x),\quad i=1,2,
\end{equation}
for all $E\in\mathcal{F}$.
\end{mytheorem}
A proof of this theorem can be found e.g in \cite{Billingsley1986}.\\\\
In the following we use $X$ for the generic variable and $x$ for a specific value of $X$. If $H_i$, $i=1,2$ is the hypothesis that X is from the statistical population with probability measure $P_i$, the mean information for discrimination in favor of $H_1$ against $H_2$ given $x\in E\in\mathcal{F}$, for $P_1$ is given by the Kullback–Leibler divergence. 
\begin{mydef}
The Kullback–Leibler (KL) divergence is given by
\begin{equation}
\begin{split}
 \mathcal{I}(P_1:P_2;E):=&\dfrac{1}{P_1(E)}\int_E\log\dfrac{f_1(x)}{f_2(x)}\mathrm{d}P_1(x)\\
 =&\left\{ \begin{matrix}  \dfrac{1}{P_1(E)}\mathlarger{\int}_Ef_1(x)\log\dfrac{f_1(x)}{f_2(x)}\mathrm{d}\lambda(x),& \mathrm{if}\ P_1(E)>0,\\
 0,& \mathrm{if}\ P_1(E)=0,   \end{matrix} \right.
\end{split}
\end{equation}
with
\begin{equation}
 \mathrm{d}P_1(x)=f_1(x)\mathrm{d}\lambda(x).
\end{equation}
When $E$ is the entire sample space $\Omega$, we shorten the notation to
$\mathcal{I}(P_1:P_2)$. For discrete sets $E$
the integral is substituted by a sum.
\label{KullbackLeiblerdivergence}
\end{mydef}
For details we refer to \cite{Kullback1997}.\\\\
Now assume that the sample space $\Omega$ is split into two disjoint sets
$E_1$ and $E_2$, $\Omega=E_1\cup E_2$. We define a statistical test procedure
to choose between hypotheses $H_1$ and $H_2$ by accepting $H_1$ if $x\in E_1$
and accepting $H_2$ if $x\in E_2$. Assuming that one of the hypotheses has to
be true we treat $H_2$ as the null hypothesis and call $E_1$ the critical
region. The following wrong test decisions can occur. 
\begin{mydef}
Incorrectly accepting $H_1$ although $H_2$ is true is
called the type I error. The probability that this error occurs is given by 
\begin{equation}
 \alpha=\Prob(x\in E_1|H_2)=P_2(E_1).
\end{equation}
\end{mydef}

\begin{mydef}
Incorrectly accepting $H_2$ although $H_1$ is true is called the type II error. The probability that this error occurs is given by
\begin{equation}
 \beta=\Prob(x\in E_2|H_1)=P_1(E_2).
\end{equation}
\end{mydef}
We assume that the test is repeated $n$-times and denote by $\mathcal{O}_n$ a
sample of $n$ independent observations. $\mathcal{O}_1$ represents a sample of
a single observation. $\beta_n$ is defined as the corresponding 
probability of an error of type II which depends on the number of
independent observations and the splitting of the probability space $\Omega$
into disjoint sets $E_1$ and $E_2$.\\\\
The following theorem demonstrates an asymptotic relation between the
KL-divergence and the minimum possible probability $\beta_n^*$ of an error of
type II with respect to all possible splittings $E_1\cup
  E_2=\Omega$ with given $\alpha=\Prob(x\in E_1|H_2)=P_2(E_1)$ \cite{Chernoff1956}. 
 \begin{mytheorem}
  For any value of $\alpha$, say $\alpha_0$, $0<\alpha_0<1$,
  \begin{equation}
   \lim_{n\rightarrow\infty}\begin{pmatrix}\beta_n^*\end{pmatrix}^{1/n}=e^{-\mathcal{I}(P_2:P_1,\mathcal{O}_1)}
  \end{equation}
 \end{mytheorem}
A proof of this theorem is given in \cite{Chernoff1956,Kullback1997}.\\\\

Assuming probability models for the outcome of a data measurement experiment
depending on experimental design parameters $\xi \in \Xi \subset
\mathbb{R}^d$,
this theorem justifies the KL-divergence to be an appropriate objective
functional for model-based computation of an optimal experimental design for
discrimination between model hypotheses. For a design with the largest
possible value of $\mathcal{I}$ the asymptotic probability $\beta_n^*$ of encountering
an error of type II  becomes minimal
with respect to all possible splittings $E_1\cup E_2=\Omega$ with given $\alpha_0$. We indicate the
dependency of the KL-divergence on the design by
$\mathcal{I}(P_2:P_1,\mathcal{O}_1;\xi)$. Our aim is to derive an algorithm to
calculate the optimal design $\hat{\xi}\in\Xi$ such that  
\begin{equation}
 \hat{\xi}=arg\max_{\xi\in\Xi}\mathcal{I}(P_2:P_1,\mathcal{O}_1;\xi).
\end{equation}
An extension of the case to test a simple null hypothesis against a simple
alternative hypothesis to the more general case of both hypotheses being
composite is generally of interest. This includes the situation to test whether given
measurement data can be explained best by one out of a finite set of \corrFIVEDOM{probability models based on} measures $P_{r_1}$, $r_1\in\{1,...,M_1\}$ parametrized by parameters $\theta_{r_1}\in\Theta_{r_1}$ where $\Theta_{r_1}
\subset \mathbb{R}^{p_{r_1}}$ is the set of all possible parameter values to parametrize $P_{r_1}$ and $M_1$ is the cardinality of the set of probability \corrFIVEDOM{models}, against the hypothesis that the measurement can best be explained by another one out of a second finite set of probability \corrFIVEDOM{models based on} measures $P_{r_2}$, parametrized
by parameters $\theta_{r_2}\in\Theta_{r_2}$ where $\Theta_{r_2}
\subset \mathbb{R}^{p_{r_2}}$ and $r_2\in\{1,...,M_2\}$.\\
By calculating
\begin{equation}\label{optexpdesign}
 \hat{\xi}=arg\max_{\xi\in\Xi}\min_{\begin{subarray}{c} r_1\in\{1,...,M_1\} \\ r_2\in\{1,...,M_2\} \end{subarray}
}\min_{\begin{subarray}{c} \theta_{r_1}\in\Theta_{r_1} \\ \theta_{r_2}\in\Theta_{r_2}\end{subarray}}\mathcal{I}(P_{r_2}(\theta_{r_2}):P_{r_1}(\theta_{r_1}),\mathcal{O}_1;\xi)
\end{equation}we can get a robust worst case estimate of an optimally discriminating design
for the case of composite null and alternative hypothesis.\\
In practical \corrFIVEDOM{applications} a simple strategy to sort different probability models into null and alternative hypothesis would be to first rank all models according to the existing measurements and then set the best model as null hypothesis and the others as alternative hypothesis. \corrFIVEDOM{The} development of a suitable and efficient strategy is subject to further work.\\
It should be noted that the presented criteria is not symmetric, i.e. the null hypothesis is favored. In the case that both hypotheses might be equally reasonable we suggest to use the symmetrized version of the KL-divergence, i.e.
\begin{equation}
 \mathcal{I}_{\mathrm{sym}}(P_2:P_1,\mathcal{O}_1;\xi)=\frac{\mathcal{I}(P_2:P_1,\mathcal{O}_1;\xi)+\mathcal{I}(P_1:P_2,\mathcal{O}_1;\xi)}{2},
\end{equation}
as optimization \corrFIVEDOM{objective} instead.
\subsection{Derivation of the optimal experimental design criterion}
\label{biodesign}
In this section we derive a numerically computable optimization objective
functional based on the framework of KL divergence. The derivation is
motivated by the requirements of biological in vitro time series experiments
modeled by kinetic ODE systems. In most situations such experiments are time
and cost consuming. Therefore a central issue is to get the most information
out of a single time series data measurement experiment taking place within a given fixed time span $[0,T_{\mathrm{end}}]$. This means that in an optimal experimental
design the most informative measurement time points for one measurement run
have to be calculated in such a way that only one measurement at one time
point can be performed. Often, an experiment cannot produce measurements in
a time continuous way. Therefore we assume that there has to be a minimal time
span $\Delta T$ for the separation of \corrFIVEDOM{subsequent} measurement time points. This contrasts to the usual approach to associate weights to a discrete or continuous time design scheme, see e.g. \cite{ATKINSON1975}.
Additionally, the initial species concentrations of the participating species
should be chosen in a most discriminating way. \\ 
A commonly used \corrFIVEDOM{experimental} practice is to combine kinetic time series measurements with
perturbation stimuli like external adding of species quantities. From the
model discrimination point of view the optimal time point of \corrFIVEDOM{perturbations} and
the 
optimal species quantities to be added should be determined. We further assume
that a measurement cannot be done at the same time as a perturbation.\\ 
In the following we translate these experimental conditions into a statistical
model. Given the measurement time-vector $t\in\mathbb{R}^n_+$ with entries
$t_i$ for the $n$ measurement time points $t_i, i\in\{1,...,n\}$ such that $t_{i+1}\geq t_i$, the ``internal'' model response vectors $y_{r_{j}}^i$ at measurement time $t_i$ for the parametrized probability measures $P_{r_j}$, $r_j\in\{1,...,M_j\}$, $j\in\{1,2$\} are given by
\begin{equation}
y_{r_{j}}^i := y_{r_{j}}(t_{i-1},t_i,y_{r_{j}}^{i-1}+c_{i-1},\theta_{r_{j}}), 
\end{equation} where $y_{r_{j}}(t_{i-1},t_i,y_{r_{j}}^{i-1}+c_{i-1},\theta_{r_{j}}) \in \mathbb{R}^{m_{r_{j}}}$ are the solutions of the initial value problems
\begin{equation}
 \dfrac{\mathrm{d}y_{r_{j}}}{\mathrm{d}\mathfrak{t}}=f_{r_{j}}^{\mathrm{rhs}}(y_{r_{j}},\theta_{r_{j}}), \quad \mathfrak{t} \in [t_{i-1},t_{i}], 
 \label{initialValueProblem}
\end{equation}with initial state $y_{r_{j}}(t^{i-1}):=y_{r_{j}}^{i-1}+c_{i-1}$ at end time $t_{i}$ where $t_0:=0$ and $c_0:=0$. The vectors $c_i  \in \mathbb{R}^{m_{\mathrm{max}}}$ denote species quantities the system can be perturbed with at time points $t_i$ where $m_{\mathrm{max}}$ is the maximum dimension of the models, i.e.
\begin{equation}
 m_{\mathrm{max}}:=\max_{j\in\{1,2\}}\max_{r\in\{1,...,M_j\}}m_{r_j}.
\end{equation} $f_{r_j}^{\mathrm{rhs}}(\cdot,\cdot)$ are the right hand side functions of the ODE models. \corrFIVEDOM{ $y_I\in\mathbb{R}^{m_{\mathrm{max}}}$ denotes the initial species concentration vector of the entire experiment  which is for all models the same, i.e. $y_{r_j}^0:= y_I$}. By $m_{\mathrm{min}}$ we denote the minimal dimension of the models, i.e
\begin{equation}
 m_{\mathrm{min}}:=\min_{j\in\{1,2\}}\min_{r\in\{1,...,M_j\}}m_{r_j}.
\end{equation}
We do not assume that $m_{\mathrm{min}}=m_{\mathrm{max}}$ therefore for a model with $m_{r_j}<m_{\mathrm{max}}$ the ``redundant'' entries in $c_i$ and $y_I$ are ``ignored''.\\\\ To each ODE model we associate an observable function $O_{r_j}:\mathbb{R}^{m_{r_j}}\rightarrow\mathbb{R}^o$ which describes an experimental observation explained by that model where $o$ denotes the dimension of the experimental observation. The expected observation $O^i_{r_j}$ of the $r_j$-th model at time point $i$ is given by
\begin{equation}
 O^i_{r_j}:=O_{r_j}(y^i_{r_j}).
\end{equation}Let \corrFIVEDOM{$O^{t_i}$} denote an observation at measurement time point $t_i$. By assuming that the measurements at successive time points $t_i$ are independent with normally distributed error vectors $\epsilon^i_{{r_j}} \in \mathbb{R}^{m}$ \corrFIVEDOM{with} zero mean
and variance functions $v_{{r_j}}(O_{r_j}^i,t_i,\theta_{{r_j}})^2$ we get for the regression models 
\corrFIVEDOM{
\begin{equation}
 O^{t_i}=O_{{r_j}}^i+\epsilon^i_{{r_j}}
\end{equation}}
the model probability densities $f_{r_j}(\cdot;\cdot)$ at measurement time point $t_i$ given by
\corrFIVEDOM{\begin{equation}
f_{r_j}(O_{t_i};O_{r_j}^{i})=\frac{1}{\sqrt{2\pi}\vert v_{i,r_j}\vert}e^{-\frac{1}{2}(O^i_{r_j}-O^{t_i})^TV_{r_j}^i(O^i_{r_j}-O^{t_i})}, 
\label{NORMALDISTRIBUTION}
\end{equation}}
with $\vert v_{i,r_j}\vert:=\prod_{k=1}^{\corrFIVEDOM{o}}v^k_{r_j}(O_{r_j}^i,t_i,\theta_{r_j})$, where $v_{r_j}^k(O_{r_j}^i,t_i,\theta_{r_j})$ denotes the $k$-th entry of the square root of the variance functions $v_{r_j}(O_{r_j}^i,t_i,\theta_{r_j})^2$, and
diagonal matrices $V^i_{r_j} \in \mathbb{R}^{o \times o}$ with diagonal entries
$\begin{bmatrix}V_{r_j}^{i}\end{bmatrix}_{kk}:=(1/v^k_{r_j}(O_{r_j}^i,t_i,\theta_{r_j}))^2$.\\\\
We generally allow different error models. The error models might dependent on the expected observations $O^i_{\corrFIVEDOM{r_j}}$, the time $t_i$ and possibly on parameters $\theta_{r_j}$.\\\\
For the sake of notational simplicity we define
\begin{equation}
 f_{r_j}(\corrFIVEDOM{O^{t_i}}):=f_{r_j}(\corrFIVEDOM{O^{t_i}};O_{r_j}^{i}).
\end{equation}
For \corrFIVEDOM{a} full measurement run containing $n$ measurement time points we get the probability density models
\begin{equation}
 f_{r_j}(O):=\prod_{i=1}^{n}f_{r_j}(\corrFIVEDOM{O^{t_i}}).
\end{equation}However, by assuming such a model probability distribution we still allow that
two measurements are separated by a time span less than $\Delta T$.\\\\
To overcome this problem we extend the probability spaces $\Omega_i=\mathbb{R}^o$ of a measurement at one measurement time point by  one-element-containing sets $\mathcal{N}_i$ to
\begin{equation}
 \tilde{\Omega}_i=\Omega_i \mathaccent\cdot\cup  \mathcal{N}_i
\end{equation}
where $\tilde{\Omega}_i$ is the disjoint union of $\Omega_i $ and $\mathcal{N}_i$.
The element of the set $\mathcal{N}_i$ with measure $P(\mathcal{N}_i)\in[0,1]$ represents the event ``no measurement'', i.e. $\widetilde{O}_{t_i} \in \mathcal{N}_i \Leftrightarrow$ ``no measurement performed at time point $t_i$''.\\\\
In order to derive measures on $\tilde{\Omega_i}$, $i=1,...,n$ that allow for a density function representation according to the Radon-Nikodym theorem
  (Theorem \ref{theorem:radon_nikodym}), we introduce the Heaviside-function
\begin{equation}
 \corrFIVEDOM{\mathcal{H}}:\mathbb{R}_+\longrightarrow [0,1]
\end{equation}with
\begin{equation}
 \corrFIVEDOM{\mathcal{H}}(t_i) = \left\{ \begin{array}{ll}
         1 & \mbox{if $t_i -t_{i-1} \geq \Delta T$}\\
         0 & \mbox{if $t_i -t_{i-1} < \Delta T$}\end{array}  \right.   
\end{equation}
By use of this Heaviside-function and $\sigma$-algebras
  $\mathcal{F}_i$, where $\mathcal{F}_i$ contains the Lebesgue measurable sets on $\Omega_i$ and additionally the union of these with the set $\mathcal{N}_i$, we define probability spaces
  $(\tilde{\Omega}_i,\mathcal{F}_i,\tilde{P}_{i,r_j})$ with measures
\begin{equation}
 \tilde{P}_{i,r_j}:E_i \in \mathcal{F}_i \mapsto \tilde{P}_{i,r_j}(E_i) \in [0,1].
\end{equation}Three cases have to be distinguished: 1.
 $E_i \subset\Omega_i$, 2. $E_i \subset\mathcal{N}_i$, 3. $E_i
 \cap\Omega_i\neq\emptyset$ and $E_i \cap\mathcal{N}_i\neq\emptyset$.\\
For case one with $E_i \subset\Omega_i$ we set 
\begin{equation}
\tilde{P}_{i,r_j}(E_i):=\corrFIVEDOM{\mathcal{H}}(t_i)\int_{E_i}f_{r_j}(\corrFIVEDOM{O^{t_i}})\mathrm{d}\corrFIVEDOM{O^{t_i}}.
\end{equation}
For case two with $E_i \subset\mathcal{N}_i$ we set
\begin{equation}
 \tilde{P}_{i,r_j}(E_i):=1-\corrFIVEDOM{\mathcal{H}(t_i)}.
\end{equation}
For case three with $E_i\cap\Omega_i\neq\emptyset$ and $E_i \cap\mathcal{N}_i\neq\emptyset$ we set
\begin{equation}
\tilde{P}_{i,r_j}(E_i):=\corrFIVEDOM{\mathcal{H}(t_i)}\int_{E_i\cap\Omega_i}f_{r_j}(\corrFIVEDOM{O^{t_i}})\mathrm{d}\corrFIVEDOM{O^{t_i}}+\left(1-\corrFIVEDOM{\mathcal{H}}(t_i)\right).
\end{equation}By introducing these modifications the probability \corrFIVEDOM{models based on measures} $\tilde{P}_{i,r_j}$ do not depend on measurements which are performed in less than $\Delta T$ time after the previous measurement any more.\\\\
To take into account that a species concentration perturbation can only be applied if no measurement is done at the same time, the same procedure is
repeated with the additional Heaviside-function 
\begin{equation}
 \corrFIVEDOM{\widetilde{\mathcal{H}}} (c_i) = \left\{ \begin{array}{ll}
         0 & \mbox{if $c_i >  0$}\\
         1 & \mbox{if $c_i = 0$.}\end{array}\right.   
\end{equation}
The measures $\tilde{P}_{i,{r_j}}$ are defined in the same way as above replacing  
\corrFIVEDOM{$\mathcal{H}(t_i)$} \corrFIVEDOM{by} \corrFIVEDOM{$\mathcal{H} (t_i)\widetilde{\mathcal{H}} (c_{i})$}.\\\\
For specific $r_1$ and $r_2$ inserting the  probability models $\tilde{P}_{i,{r_1}}$ and respectively $\tilde{P}_{i,{r_2}}$  into the KL-divergence (Definition \ref{KullbackLeiblerdivergence}) using $\lambda_i:=\tilde{P}_{i,r_2}$ and the additivity of the KL divergence for independent events one gets the following expression
\corrFIVEDOM{\begin{equation}\begin{split} 
 \mathcal{I}(P_{r_2}:P_{r_1},\mathcal{O}_1)=\sum_{i=1}^n& \left[\int\mathcal{H} (t_i)\widetilde{\mathcal{H}} (c_{i})f_{r_2}(O^{t_i})
\log\begin{Bmatrix}\dfrac{\mathcal{H} (t_i)\widetilde{\mathcal{H}} (c_{i})f_{r_2}(O^{t_i})}{\mathcal{H} (t_i)\widetilde{\mathcal{H}} (c_{i})f_{r_1}(O^{t_i})}\end{Bmatrix} \mathrm{d}O^{t_i}\right. + \\
&\qquad\qquad\qquad\quad\left.\left( 1-\mathcal{H} (t_i)\widetilde{\mathcal{H}}(c_{i}) \right)\cdot\log{\dfrac{ (1-\mathcal{H} (t_i)\widetilde{\mathcal{H}}(c_{i}))}{ (1-\mathcal{H} (t_i)\widetilde{\mathcal{H}}(c_{i}))}} \right], 
\end{split}
\end{equation}}
where $c_n:=0$.
With $\log(1)=0$ this simplifies to
\corrFIVEDOM{\begin{equation}\begin{split}
 \mathcal{I}(P_{r_2}:P_{r_1},\mathcal{O}_1)=\sum_{i=1}^n\mathcal{H} (t_i)\widetilde{\mathcal{H}} (c_{i})\int f_{r_2}(O^{t_i})\cdot
\log\begin{Bmatrix}\dfrac{f_{r_2}(O^{t_i})}{f_{r_1}(O^{t_i})}\end{Bmatrix} \mathrm{d}O^{t_i}.\end{split}
\label{NORMALDISTRIBUTION2}
\end{equation}}
By inserting the normal distribution (\ref{NORMALDISTRIBUTION}) in (\ref{NORMALDISTRIBUTION2}) one gets
\corrFIVEDOM{\begin{equation}
\begin{split}
 \mathcal{I}(P_{r_2}:P_{r_1},\mathcal{O}_1&)=\sum_{i=1}^n\mathcal{H} (t_i)\widetilde{\mathcal{H}} (c_{i})\cdot\\
 &\int f_{r_2}(O^{t_i})\cdot\log\begin{Bmatrix}\dfrac{\frac{1}{\sqrt{2\pi}\vert v_{i,r_2}\vert}e^{-\frac{1}{2}(O^i_{r_2}-O^{t_i})^TV_{r_2}^i(O^i_{r_2}-O^{t_i})}}{\frac{1}{\sqrt{2\pi}\vert v_{i,r_1}\vert}e^{-\frac{1}{2}(O^i_{r_1}-O^{t_i})^TV_{r_1}^i(O^i_{r_1}-O^{t_i})}}\end{Bmatrix} \mathrm{d}O^{t_i}.
\end{split}\end{equation}}
This is equivalent to
\corrFIVEDOM{\begin{equation}
\mathcal{I}(P_{r_2}:P_{r_1},\mathcal{O}_1)=\sum_{i=1}^n\mathcal{H} (t_i)\widetilde{\mathcal{H}} (c_{i})\left(\sum_{k=1}^o \log\frac{ v^k_{r_1}(O^i_{r_1},t_i,\theta_{r_1})}{ v^k_{r_2}(O^i_{r_2},t_i,\theta_{r_2})}+A^i_k\right),
\label{fullEquation}
\end{equation}}
with 
\corrFIVEDOM{\begin{equation}
\begin{split}
 A^i_k:=&\frac{1}{2}\int \left\lfloor f_{r_2}(O^{t_i})\right\rfloor_k\cdot\\   
   &\left[-\frac{1}{(v^k_{r_2}(O^i_{r_2},t_i,\theta_{r_2}))^2}\left( \left\lfloor O^i_{r_2}\right\rfloor_k^2-2\left\lfloor O^i_{r_2}\right\rfloor_k\left\lfloor O^{t_i}\right\rfloor_k+\left\lfloor O^{t_i}\right\rfloor_k^2\right)\right.+ \\ 
   &\left. \frac{1}{(v^k_{r_1}(O^i_{r_1},t_i,\theta_{r_1}))^2}\left( \left\lfloor O^i_{r_1}\right\rfloor_k^2-2\left\lfloor O^i_{r_1}\right\rfloor_k \left\lfloor O^{t_i}\right\rfloor_k+ \left\lfloor O^{t_i}\right\rfloor_k^2\right) \right] \left\lfloor\mathrm{d}O^{t_i}\right\rfloor_k.
\end{split}
\end{equation}}
and where $\left\lfloor O \right\rfloor_k$ gives the $k-th$ entry of the observation vector $O$. $A^i_k$ reduces using the well known moments of the normal distribution to
\begin{equation}
\begin{split}
 &A^i_k=\frac{1}{2}\left[-\frac{1}{(v^k_{r_2}(O^i_{r_2},t_i,\theta_{r_2}))^2}\left( \left\lfloor O^i_{r_2}\right\rfloor_k^2-2\left\lfloor O^i_{r_2}\right\rfloor ^2_k+\left\lfloor O^i_{r_2}\right\rfloor^2_k+(v^k_{r_2}(O^i_{r_2},t_i,\theta_{2}))^2\right)\right. \\ &\left. +\frac{1}{(v^k_{r_1}(O^i_{r_1},t_i,\theta_{r_1}))^2}\left( \left\lfloor O^i_{r_1}\right\rfloor_k^2-2\left\lfloor O^i_{r_1}\right\rfloor_k\left\lfloor O^i_{r_2}\right\rfloor_k+\left\lfloor O^i_{r_2}\right\rfloor^2_k+(v^k_{r_2}(O^i_{r_2},t_i,\theta_{r_2}))^2\right) \right].
\end{split}
\end{equation}
This further simplifies to
\begin{equation}
\begin{split}
 A^i_k=\frac{1}{2}\left[-1  +\frac{\left( \left\lfloor O^i_{r_1}\right\rfloor_k^2-2\left\lfloor O^i_{r_1}\right\rfloor_k\left\lfloor O^i_{r_2}\right\rfloor_k+\left\lfloor O^i_{r_2}\right\rfloor^2_k+(v^k_{r_2}(O^i_{r_2},t_i,\theta_{r_2}))^2\right)}{(v^k_{r_1}(O^i_{r_1},t_i,\theta_{r_1}))^2} \right].
\end{split}
\end{equation}
Substituting $A^i_k$ back into (\ref{fullEquation}) we get
\corrFIVEDOM{\begin{equation}
\begin{split}
  \mathcal{I}(P_{r_2}:P_{r_1},\mathcal{O}_1)&=\sum_{i=1}^n\mathcal{H} (t_i)\widetilde{\mathcal{H}} (c_{i})\sum_{k=1}^o\left( \log{\left(\frac{v^k_{r_1}(O^i_{r_1},t_i,\theta_{r_2})}{v^k_{r_2}(O^i_{r_2},t_i,\theta_{r_1})}\right)}+\right.\\
  &\left.\frac{1}{2}\left[-1  +\frac{ \left(\left\lfloor O^i_{r_1}\right\rfloor_k-\left\lfloor O^i_{r_2}\right\rfloor_k\right)^2+(v^k_{r_2}(O^i_{r_2},t_i,\theta_{r_2}))^2}{(v^k_{r_1}(O^i_{r_1},t_i,\theta_{r_1}))^2} \right]\right).
\end{split}
\end{equation}}
This reduces to
\corrFIVEDOM{\begin{equation}
  \begin{split}
 \mathcal{I}(P_{r_2}:P_{r_1},\mathcal{O}_1)=&\dfrac{1}{2}\sum_{i=1}^n  \mathcal{H} (t_i)\widetilde{\mathcal{H}} (c_{i})\cdot\\
 &\left( \sum_{k=1}^o \left[ \dfrac{(v^k_{r_2}(O^i_{r_2},t_i,\theta_{r_2}))^2+\begin{pmatrix}\left\lfloor O_{r_2}^{i}\right\rfloor_k-\left\lfloor O_{r_1}^{i}\right\rfloor_k\end{pmatrix}^2}{(v^k_{r_1}(O^i_{r_1},t_i,\theta_{r_1}))^2} \right.\right.\\ &\qquad\qquad\qquad\qquad \left.\left. -2\log\begin{pmatrix}\dfrac{v^k_{r_2}(O^i_{r_2},t_i,\theta_{r_2})}{v^k_{r_1}(O^i_{r_1},t_i,\theta_{r_1})}\end{pmatrix}\right]-o \right) .
 \end{split}
 \label{optimizationCriterion}
\end{equation}}
This criterion has to be maximized with respect to the initial concentration
\corrFIVEDOM{vector} $y_I$, the measurement time point \corrFIVEDOM{vector} $t$ and the system perturbation \corrFIVEDOM{vector} $c$, thus $\xi:=(y_I,t,c)\in\Xi\subset\mathbb{R}^d$.\\ For
our optimal experimental design we generally start with a large number of
measurement time points. By use of the \corrFIVEDOM{Heaviside} functions the number of measurement
time points \corrFIVEDOM{is} reduced such that for $t_i-t_{i-1}<\Delta T$ the
corresponding measurement time point is ``turned off''.\\\\ \corrFIVEDOM{These
  Heaviside-functions $\mathcal{H} (\cdot)$ and $\widetilde{\mathcal{H}}
  (\cdot)$ can be replaced by any appropriate continuously differentiable
  switching functions with range space $[0,1]$.}\\
It should be noted that we assume that we have the same time discretization
for measurements of different species and the addition of further species
\corrFIVEDOM{quantities}. This assumption is practical
\corrFIVEDOM{especially} for the application to in vitro experiments performed
by \corrFIVEDOM{biologists}. \corrFIVEDOM{For introducing arbitrary generic
  controls we need a more general formulation of time schemes,
  i.e. simultaneously time schemes which are independent of each other. One is
  associated with the controls, others may be associated with distinct
  observables which might be measured indepently. The incorporation to the
  presented framework is subject of further work.}
\section{Solution of the max-min optimization problem}
\label{maxminTheory}
\begin{figure}[ht]
 \begin{center}
   \includegraphics[width=0.7\columnwidth]{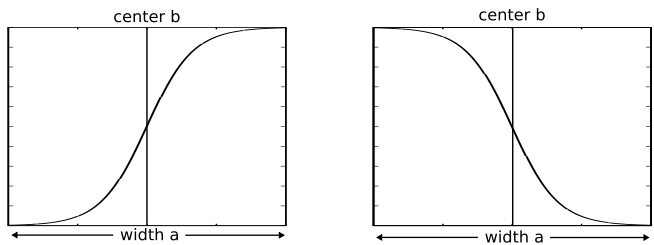}
    \end{center}
 \caption{Switching functions: the left switching function is used to guarantee that only one measurement is done at a time point, the right one is used to guarantee that if a perturbation is done at a time point no measurement is done at the same time point.}
 \label{tangh}
  \end{figure}

We formally state now the experimental design optimization problem $\mathbf{P}_\Theta$:
\begin{equation}
 \max_{(\tau,\xi)\in\Xi\subset\mathbb{R}^{d+1}}\tau
 \label{naturalMINMAXproblem}
\end{equation}
subject to 
\begin{equation*}
\renewcommand*{\arraystretch}{1.5}
 \begin{array}{c}
   \min\limits_{\begin{subarray}{c} \theta_{r_1}\in\Theta_{r_1} \\ \theta_{r_2}\in\Theta_{r_2} \end{subarray}} \mathcal{I}(P_{r_2}(\theta_{r_2}):P_{r_1}(\theta_{r_1}),\mathcal{O}_1;\xi)-\tau\geq0, \ r_1\in\{1,...,M_1\},\ r_2\in\{1,...,M_2\},\\
   \sum\limits_{i=1}^n \Delta t_i = T_{\mathrm{end}},\\
   y_{I}^{\mathrm{min}} \leq y_I \leq y_{I}^{\mathrm{max}},\\
   0\leq \Delta  t\leq t^{\mathrm{max}},\\
   0\leq c\leq c^{\mathrm{max}},
\end{array}
\end{equation*}
with $\xi:=(y_I,\Delta t,c)\in\mathbb{R}^d$, $\Delta t_i:=t_i-t_{i-1}$ and $\Theta:=\{(\Theta_{r_1},\Theta_{r_2})|r_j\in\{1,...,M_j\}\}$. The auxiliary variable $\tau\in\mathbb{R}$ is used to transform the maxmin optimization problem (\ref{optexpdesign}) \corrFIVEDOM{to} a maximization problem with an infinite number of \corrFIVEDOM{inequality} constraints. The remaining constraints model the feasible range of experimental
setups. \\\\
To solve optimization problem (\ref{naturalMINMAXproblem}) numerically by applying \corrFIVEDOM{efficient} derivative based algorithms we replace the Heaviside functions \corrFIVEDOM{$\mathcal{H}(t_i)$} and \corrFIVEDOM{$\widetilde{\mathcal{H}}(c_i)$} in (\ref{optimizationCriterion}) by continuously differentiable approximations, parametrized hyperbolic tangent functions of the form 
\begin{equation}
 \corrFIVEDOM{\mathcal{H}'}(t_i)=\frac{\tanh(\frac{6(\Delta t_i-b_1)}{a_1})+1}{2}\quad \mathrm{and} \quad \corrFIVEDOM{\widetilde{\mathcal{H}}'}(c_i)=\frac{\tanh(-\frac{6( c_i-b_2)}{a_2})+1}{2}.
\end{equation}The parameters $a_{1,2}$ characterize the width of the transition region between $0$ and $1$. The parameters $b_{1,2}$ determine the center of the transition region (see Figure \ref{tangh}). By setting the parameters in an adequate way arbitrarily close approximations of the Heaviside functions can be generated.
A different approach to handle the discontinuous Heaviside functions \corrFIVEDOM{would be} to introduce \corrFIVEDOM{binary} variables and treat the resulting problem as  Mixed Integer Nonlinear Programming problem. The drawback of this approach \corrFIVEDOM{is} that \corrFIVEDOM{its} solution \corrFIVEDOM{can become} very expensive. There seems \corrFIVEDOM{to be little} theoretical work in literature \corrFIVEDOM{on} Mixed Integer \corrFIVEDOM{maxmin} problems and \corrFIVEDOM{an efficient} solution strategy \corrFIVEDOM{is not obvious} in that case.
\\
In literature\corrFIVEDOM{,} problems as (\ref{naturalMINMAXproblem}) fall into the class of semi-infinite inequality and equality constrained optimization problems (\textbf{SIECP}) \cite{Polak1997}.\\ 
 Several methods to solve such \textbf{SIECP} problems are available, an overview can be found in \cite{Hettich1993,Polak1997}. We choose the method of outer approximation \cite{Salmon1968,Shimizu1980,Polak1997}, whose origin can be traced back to cutting plane methods for convex problems \cite{Polak1997}. This approach is beneficial in the presence of a complex inner problem, in our case the robustification against the parameters $\theta_{r_1} \mathrm{and}\ \theta_{r_2}$. The outer approximation algorithm \corrFIVEDOM{iteratively} solves discretized finite counterparts $\mathbf{P}_{\widetilde{\Theta}^N}$ of the semi-infinite problem $\mathbf{P}_{\Theta}$ in each step $N$\corrFIVEDOM{,} successively refining the discretization $\widetilde{\Theta}^N$ until a sufficient approximation of the original problem $\mathbf{P}_{\Theta}$ is reached. For problem (\ref{naturalMINMAXproblem}) this means that in each iteration of the outer approximation scheme $\Theta_{r_1}$ and $\Theta_{r_2}$ \corrFIVEDOM{are} replaced by finite approximations $\widetilde{\Theta}^N_{r_1}$ and respectively $\widetilde{\Theta}^N_{r_2}$ with $\widetilde{\Theta}^N:=\{(\widetilde{\Theta}^N_{r_1},\widetilde{\Theta}^N_{r_2})|r_j\in\{1,...,M_j\}\}$. This relation between the semi-infinite problem and an infinite sequence of finite problems can be formalized in the theory of consistent approximations and epi-convergence \cite{Polak1969,Polak1987,Polak1993,Polak1997}. \\\\
We use a modified version of Algorithm 3.6.4 in \cite{Polak1997}  where ``Step 1.''\corrFIVEDOM{,} the calculation of augmenting vectors $\tilde{\theta}^{N+1}_{r_1}$ and $\tilde{\theta}^{N+1}_{r_1}$ to construct
\begin{equation}
 \widetilde{\Theta}^{N+1}_{r_1}:=\{\tilde{\theta}^{N+1}_{r_1}\}\cup\widetilde{\Theta}^N_{r_1}\quad \mathrm{and}\quad\widetilde{\Theta}^{N+1}_{r_2}:=\{\tilde{\theta}^{N+1}_{r_2}\}\cup\widetilde{\Theta}^N_{r_2},
\label{StepOneFormula1}
\end{equation}
is \corrFIVEDOM{realized} by 
\begin{equation}
 (\tilde{\theta}^{N+1}_{r_1},\tilde{\theta}^{N+1}_{r_2}):=arg\min_{\begin{subarray}{c}\theta_{r_1}\in\Theta_{r_2} \\ \theta_{r_2}\in\Theta_{r_2}\end{subarray}}\mathcal{I}(P_{r_2}(\theta_{r_2}):P_{r_1}(\theta_{r_1}),\mathcal{O}_1;\xi_N),
\label{StepOneFormula2}
\end{equation}
with $\xi_N$ denoting \corrFIVEDOM{a} \corrFIVEDOM{locally} optimal design of the \corrFIVEDOM{previous} problem $\mathbf{P}_{\widetilde{\Theta}^N}$.
The algorithmic scheme is as follows:
\begin{myalgorithm}{}\mbox{}\\
\rule{\textwidth}{.3pt}
\begin{itemize}
 \item{Data.} Choose $\xi_0\in\Xi$ and a sequence $\{\epsilon_N\}^\infty_{N=1}$ with $\epsilon_N>0$ and $\epsilon_N\downarrow 0$.
\item{Step 0.} Set $N=1$, set $\widetilde{\Theta}_{r_j}^0:=\emptyset$.
\item{Step 1.} Calculate $\widetilde{\Theta}^N$ according to (\ref{StepOneFormula1}) and (\ref{StepOneFormula2}).  
\item{Step 2.} Calculate approximate solution of $\mathbf{P}_{\widetilde{\Theta}^N}$ such that
\begin{equation}
\corrFIVEDOM{\Psi}_{\widetilde{\Theta}_N}((\tau_N,\xi_N))\geq -\epsilon_N, 
\end{equation}
and the \corrFIVEDOM{equality and inequality} constraints in problem \corrFIVEDOM{$\mathbf{P}_{\widetilde{\Theta}^N}$} are fulfilled up to $\epsilon_N$.
\item{Step 3.} Replace $N$ by $N+1$, and goto Step 1.
\end{itemize}
\rule{\textwidth}{.3pt}
\label{relaxationAlgorithm}
\end{myalgorithm}
\noindent$\corrFIVEDOM{\Psi}_{\widetilde{\Theta}_N}(\cdot):\mathbb{R}^{d+1}\rightarrow\mathbb{R}_{\leq0}$ denotes the optimality function associated to problem $\mathbf{P}_{\widetilde{\Theta}^N}$, see Theorem 2.2.24 in \cite{Polak1997}. The optimality function $\corrFIVEDOM{\Psi}_{\widetilde{\Theta}_N}(\cdot)$ is always non positive and directly related to the first order generalized Karush-Kuhn-Tucker (KKT)
conditions, i.e. $\corrFIVEDOM{\Psi}_{\widetilde{\Theta}_N}((\tau,\xi))=0$ if evaluated at a generalized  KKT point, see Theorem 2.2.19 in \cite{Polak1997}.\\
Assuming that $\mathcal{I}(P_{r_2}(\theta_{r_2}):P_{r_1}(\theta_{r_1}),\mathcal{O}_1;\xi)$ and $\nabla_\xi \mathcal{I}(P_{r_2}(\theta_{r_2}):P_{r_1}(\theta_{r_1}),\mathcal{O}_1;\xi)$ are Lipschitz continuous on bounded sets \corrFIVEDOM{with} respect \corrFIVEDOM{to} $\xi$ and $\theta_{r_j}$ and $\Theta_{r_j}$ are compact any accumulation point of Algorithm \ref{relaxationAlgorithm} fulfills the generalized KKT conditions\corrFIVEDOM{,} \corrFIVEDOM{c}ompare Theorem 3.6.5 in \cite{Polak1997}.\\\\
To calculate $\tilde{\theta}^{N+1}_{r_j}$ in Step 1. of Algorithm \ref{relaxationAlgorithm} we use on heuristic base a simple random search approach coupled to a local optimization method, i.e. we have randomly generated $P$ different start values in $\Theta_{r_j}$, from which we have started the local optimization  method \corrFIVEDOM{for parameter estimation}. The best value out of the $P$ trials is chosen to augment the set $\widetilde{\Theta}^{N}_{r_j}$. Of course there are more sophisticated approaches  to search for a global minimum for a review see e.g. \cite{Arora1995}, but at this point an effective calculation of Step 1. of Algorithm \ref{relaxationAlgorithm} was not our primary goal. For the local \corrFIVEDOM{parameter optimization} we use the same optimization method as for Step 2. in Algorithm \ref{relaxationAlgorithm}.\\\\ 
In our implementation we use a fixed $\epsilon$ at the desired final accuracy to solve problem $\mathbf{P}_{\widetilde{\Theta}^N}$ in every loop of Step 2. of Algorithm \ref{relaxationAlgorithm}, i.e. $\epsilon_N=\epsilon$, $N>0$. In that way Step 1. of Algorithm  \ref{relaxationAlgorithm} gives a worst case estimate of the KL divergence $\mathcal{I}(P_{r_2}(\theta_{r_2}):P_{r_1}(\theta_{r_1}),\mathcal{O}_1;\xi_N)$ \corrFIVEDOM{for} the current design $\xi_N$ \corrFIVEDOM{with} respect \corrFIVEDOM{to} $\theta_{r_j}$. \corrFIVEDOM{Therefore for practical application the algorithm can be stopped if the worst case estimate of KL divergence for the current design $\xi_N$ is big enough although no local optima might be achieved during optimization}. As stopping criterion of Algorithm \ref{relaxationAlgorithm} we use: 
\begin{mySTOP}\label{stopCriterionLabel}\mbox{}\\Stop after Step 1. of Algorithm \ref{relaxationAlgorithm}, if 
\begin{equation}
\begin{split}
 \delta\geq &\corrFIVEDOM{\min_{\begin{subarray}\ r_1\in\{1,...,M_1\} \\ r_2\in\{1,...,M_2\}\end{subarray}}}\min_{\begin{subarray}{c}\theta_{r_1}\in\widetilde{\Theta}^{N-1}_{r_1} \\  \theta_{r_2}\in\widetilde{\Theta}^{N-1}_{r_2}\end{subarray}}\mathcal{I}(P_{r_2}(\theta_{r_2}):P_{r_1}(\theta_{r_1}),\mathcal{O}_1;\xi_{N-1})- \\
&\corrFIVEDOM{\min_{\begin{subarray}\ r_1\in\{1,...,M_1\} \\ r_2\in\{1,...,M_2\}\end{subarray}}}\min_{\begin{subarray}{c}\theta_{r_1}\in\Theta_{r_1} \\ \theta_{r_2}\in\Theta_{r_2} \end{subarray}}\mathcal{I}(P_{r_2}(\theta_{r_2}):P_{r_1}(\theta_{r_1}),\mathcal{O}_1;\xi_{N-1})=:\Delta_{RG},
\end{split}
\label{SToptCriterionNew}
\end{equation}
where $\delta$ is a small positive constant, then consider $(\tilde{\theta}^{N+1}_{r_1},\tilde{\theta}^{N+1}_{r_2})$, $r_1\in\{1,...,M_1\}$, $r_2\in\{1,...,M_2\}$
 and $\xi_{N}$ as (approximate) solutions of problem $\mathbf{P}_{\Theta}$, else goto Step 2. and calculate a new design $\xi_{N+1}$.
\end{mySTOP}
\noindent This stop criterion \corrFIVEDOM{has} also \corrFIVEDOM{been} used in \cite{Shimizu1980,Pronzato1988}. We call the  distance $\Delta_{RG}$ given by (\ref{SToptCriterionNew}),
\corrFIVEDOM{robustification} gap.
\subsection{Numerical solution of the problem $\mathbf{P}_{\widetilde{\Theta}^N}$}
\label{numericRealization}
\begin{figure}[h]
      \begin{center}
      \includegraphics[width=0.7\columnwidth]{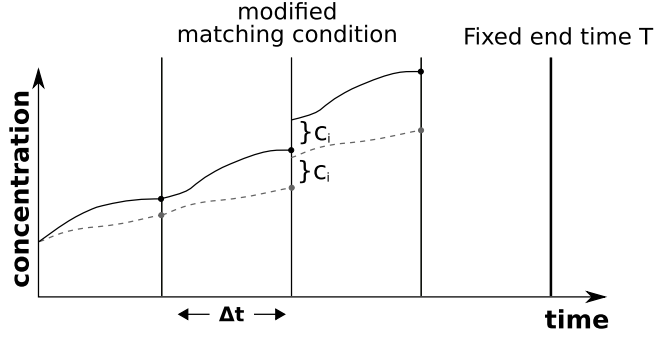}
      \end{center}
    \caption{Scheme of the multiple shooting setup for computing the experimental design \corrFIVEDOM{to discriminate two models}. \corrFIVEDOM{A} dot denotes one measurement time point. The black solid line \corrFIVEDOM{corresponds to} model 1 and the gray dashed one \corrFIVEDOM{to} model 2.}
    \label{setup}
 \end{figure}

We have implemented the resulting optimization problem in a multiple shooting setup (see for example \cite{Stoer2002,Bock1984,Bock1987}). In our multiple shooting setup the whole integration interval $[0, T_{\mathrm{end}}]$ is subdivided into several subintervals by introducing auxiliary multiple shooting node variables $s_{{r_j},i,l}$, $j\in\{1,2\}, \ r_j\in\{1,...,M_j\},\ i\in\{1,...,n\}, \ l\in\{1,...,N\}$, on each of which an independent initial value problem is solved. \corrFIVEDOM{Each} end  point of a subinterval corresponds to one measurement time point. Matching conditions which enter the optimization
problem as additional equality constraints assure continuity of the state trajectory
from one subinterval to the next.\\\\
To incorporate the perturbations $c$ matching conditions
\begin{equation}
s_{r_j,i,l}-y_{r_j}(t_{i-1},t_i,s_{r_j,i-1,l},\tilde{\theta}_{r_j}^l)=0, 
\end{equation} $s_{r_j,i,l}$, denoting the multiple shooting nodes with $s_{r_j,0,l}=y_I$ are modified to
\begin{equation}
\begin{split}
s_{r_j,i,l}-y_{r_j}(t_{i-1},t_i,s_{j,i-1,l},\tilde{\theta}_{r_j}^l)=c_i, &\quad i\in\{1,...,n-1\},\\
s_{r_j,n,l}-y_{r_j}(t_{n-1},t_n,s_{j,n-1,l},\tilde{\theta}_{r_j}^l)=0.
\end{split}
\label{modifedCconstraints}
\end{equation}
A graphical scheme of the multiple shooting setup is shown in Figure \ref{setup}.\\  Instead of evaluating the objective functional (\ref{optimizationCriterion}) by use of the values $y_{r_j}^i$, given by the solution of the initial value problem (\ref{initialValueProblem}), depending on the parameters $\tilde{\theta}_{r_j}$, (\ref{optimizationCriterion}) is evaluated by use of the auxiliary multiple shooting node variables $s_{r_j,i,l}$ replacing the values $y_{r_j}^i$ with $s_{r_j,i,l}$ respectively. The dependency of (\ref{optimizationCriterion}) on $s_{r_j,i,l}$ is indicated by $\mathcal{I}(s_{r_1,\cdot,l},s_{r_2,\cdot,l})$.
\\
 The overall optimization problem can be stated as
\begin{equation}
\max_{\tau,y_I,\Delta t, c,s}\tau
\label{naturalMINMAXproblemREFORMULATEDFULL}
\end{equation}
subject to
\begin{equation}
\renewcommand*{\arraystretch}{1.5}
 \begin{array}{c}
\frac{\mathrm{d}y_{r_j}}{\mathrm{d}\mathfrak{t}}=f^{\mathrm{rhs}}_{r_j}(y,\tilde{\theta}_{r_j}^l),  \quad \mathfrak{t}\in [t_{i-1},t_i],\quad y_{r_j}(t_{i-1}):=s_{r_j,i-1,l}\\
   s_{r_j,i,l}-y_{r_j}(t_{i-1},t_i,s_{r_j,i-1,l},\tilde{\theta}_{r_j}^l)=c_i,\\
    s_{r_j,n,l}-y_{r_j}(t_{n-1},t_n,s_{r_j,n-1,l},\tilde{\theta}_{r_j}^l)=0,\\
   s_{r_j,0,l}=y_I,\\   
   y_{I}^{\mathrm{min}} \leq  y_I \leq y_{I}^{\mathrm{max}},\\
   0\leq \Delta  t\leq t^{\mathrm{max}},\\
   0\leq c\leq c^{\mathrm{max}},\\ 
   s_{r_j,i,l}^{\mathrm{min}}\leq s_{r_j,i,l} \leq s_{r_j,i,l}^{\mathrm{max}},\\
   \sum_{i=1}^n \Delta t_i = T_{\mathrm{end}},\\
   \mathcal{I}(s_{r_1,\cdot,l},s_{r_2,\cdot,l})-\tau\geq 0, \ r_1\in\{1,...,M_1\},\ r_2\in\{1,...,M_2\},
  \end{array}
  \label{optiProblem}
\end{equation}
with $j\in\{1,2\}, \ r_j\in\{1,...,M_j\},\ i\in\{1,...,n\}, \ l\in\{1,...,N\}$.\\\\
We have  implemented this problem within the interior point optimization package IPOPT \cite{Waechter2002,Waechter2006}, using the linear solver MA27 \cite{HSL}. Usually for the solution of the KKT system within the direct multiple shooting ansatz the linear system is solved utilizing tailored structure-exploitation, e.g., condensing or high rank updates. See for example \cite{Leineweber1998}. Since speed aspects are not our primary concern we rely on the sparse solver MA27 instead of developing a tailored solver for this problem class \corrFIVEDOM{at the current stage}.\\ All derivatives up to second order, which are used for the calculations of the Hessian needed for a robust performance of IPOPT, are calculated by automatic differentiation using CppAD, \cite{Bell2008,Bell2010}.\\ For the solution of the differential equations within the optimization problem, which are commonly stiff in chemical and biochemical
applications, we have implemented a fully variable step, variable order (order $1$ to $6$), Backward differentiation formulae (BDF) method, based on Nordsiek array polynomial interpolation similar to the EPISODE BDF method by Byrne and Hindmarsh \cite{Byrne1975}, but with the step size selection strategy of Calvo and R\'andez \cite{Calvo1993}.\\
For the generation of sensitivities we have adopted the sophisticated principles of internal numerical differentiation developed by Albersmeyer and Bock \cite{Albersmeyer2008a,Albersmeyer2010} in forward and adjoint mode.\\ The idea of this principle is instead of calculating the sensitivities by use of the sensitivity differential equation, to directly differentiate the BDF integration scheme by automatic differentiation, which we implemented using CppAD \cite{Bell2008,Bell2010}.\\
According to some notes in the PhD thesis of Albersmeyer \cite{Albersmeyer2010} we also have implemented the possibility to control the step size scheme not only by the local truncation error of the nominal trajectory but as well by the local truncation error of the sensitivities generated by the forward mode of automatic differentiation with respect to the sensitivity differential equation, which has shown by numerical experience to improve the robustness of the optimization approach.\\
A different approach would be to use collocation, i.e. to \corrFIVEDOM{incorporate a full discretization of the} ODEs \corrFIVEDOM{into} the optimization problem, see e.g. \cite{Biegler2001}. Since the \corrFIVEDOM{kinetic} ODE systems \corrFIVEDOM{in the focus of our applications} are usually stiff\corrFIVEDOM{, we prefer adaptive time integration.}
\subsection{Stabilizing homotopy method for subsequent \corrFOURDOM{$\mathbf{P}_{\widetilde{\Theta}^{N+1}}$}}
\label{HomotopyStrategy}
By solving the subsequent optimization problems \corrFOURDOM{$\mathbf{P}_{\widetilde{\Theta}^{N+1}}$} with an interior point code like \corrFIVEDOM{IPOPT} \cite{Waechter2002,Waechter2006} initialized with primal and dual variables of the previous problem \corrFIVEDOM{or with primal variables only}, one often observes that the new solution may differ significantly from the previous. This is due to the fact that the solution of the previous problem $\mathbf{P}_{\widetilde{\Theta}^{N}}$ is infeasible for $\mathbf{P}_{\widetilde{\Theta}^{N+1}}$ and thus the algorithm tries to find a feasible state before it proceeds to find a new optimum. This behavior is not desired in the context of an outer approximation algorithm, because convergence of the algorithm may be slowed down \corrFIVEDOM{significantly}. This circumstance originates from a jumping between vicinities of distinct local maxima of problem (\ref{naturalMINMAXproblem}). The discretization $\widetilde{\Theta}^{N}$ of the robustification space may not be equally adequate for \corrFIVEDOM{different} local \corrFIVEDOM{maxima}. To overcome this problem we have implemented a heuristic homotopy method to gradually introduce the additional constraints  
\begin{equation}
 g_{r_1,r_2}(\tau,\xi)_{N+1}:=\mathcal{I}(s_{r_1,\cdot,N+1},s_{r_2,\cdot,N+1})-\tau\geq 0,
\end{equation}
of problem $\mathbf{P}_{\widetilde{\Theta}^{N+1}}$. We replace $g_{r_1,r_2}(\tau,\xi)_{N+1}$ by 
\begin{equation}
 \tilde{g}_{r_1,r_2}(\tau,\xi;\kappa)_{N+1}:=\mathcal{I}(s_{r_1,\cdot,N+1},s_{r_2,\cdot,N+1})-\tau+(1-\kappa)\rho\geq 0, 
\end{equation}
with homotopy parameter $\kappa\in [0,1]$ and $\rho$ is a constant which has to be set such that $\tilde{g}_{r_1,r_2}(\tau,\xi;\kappa)_{N+1}$ are inactive for $\kappa=0$ at the initial design $\xi_N$. We choose $\rho$ to be
\begin{equation}
\begin{split}
\rho:= K\cdot\max_{\begin{subarray}{c} r_1\in\{1,...,M_1\} \\  r_2\in\{1,...,M_2\} \end{subarray}} &\left( \min_{\begin{subarray}{c}\theta_{r_1}\in\widetilde{\Theta}^{N}_{r_1} \\ \theta_{r_2}\in\widetilde{\Theta}^{N}_{r_2}\end{subarray}}\mathcal{I}(P_{r_2}(\theta_{r_2}):P_{r_1}(\theta_{r_1}),\mathcal{O}_1;\xi_{N})-      \right. \\
&\ \ \ \left.   \min_{\begin{subarray}{c}\theta_{r_1}\in\Theta_{r_1} \\ \theta_{r_2}\in\Theta_{r_2} \end{subarray}}
 \mathcal{I}(P_{r_2}(\theta_{r_2}):P_{r_1}(\theta_{r_1}),\mathcal{O}_1;\xi_{N})\right).
\end{split}
\end{equation}
$K$ is a save guard factor we set \corrFIVEDOM{empirically} to $K=1.4$, which worked well in practice \corrFIVEDOM{for our examples}.
For $\kappa=0$ the augmented optimization problem should be easily solvable within a few iterations by performing a warm start from the solution of the previous problem.
By increasing the homotopy parameter to $\kappa=1$ the additional constraint is gradually introduced, which leads to a sequence of easily solvable subproblems whose solutions stay in the vicinity of the solution of the previous problem $\mathbf{P}_{\widetilde{\Theta}^{N}}$. A similar homotopy strategy can be found e.g. in \cite{Perez2009}.
\section{Numerical results}
\label{numericalResults}
We have applied the algorithm developed in Section \ref{theory} and \ref{maxminTheory} to two example problems for which we present results in the following section. For the \corrFIVEDOM{purpose} of illustration we restrict ourself to the case that each hypothesis comprise only one model whereby we assume that only the alternative hypothesis is composite. We also assume that each species is ``directly'' measurable. We treat model 1 as null hypothesis and model 2 as alternative hypothesis.
\subsection{Discriminating design for two models describing glycolytic oscillations}
\begin{figure}[ht]
\begin{center}
\subfigure{\includegraphics[width=0.49\columnwidth]{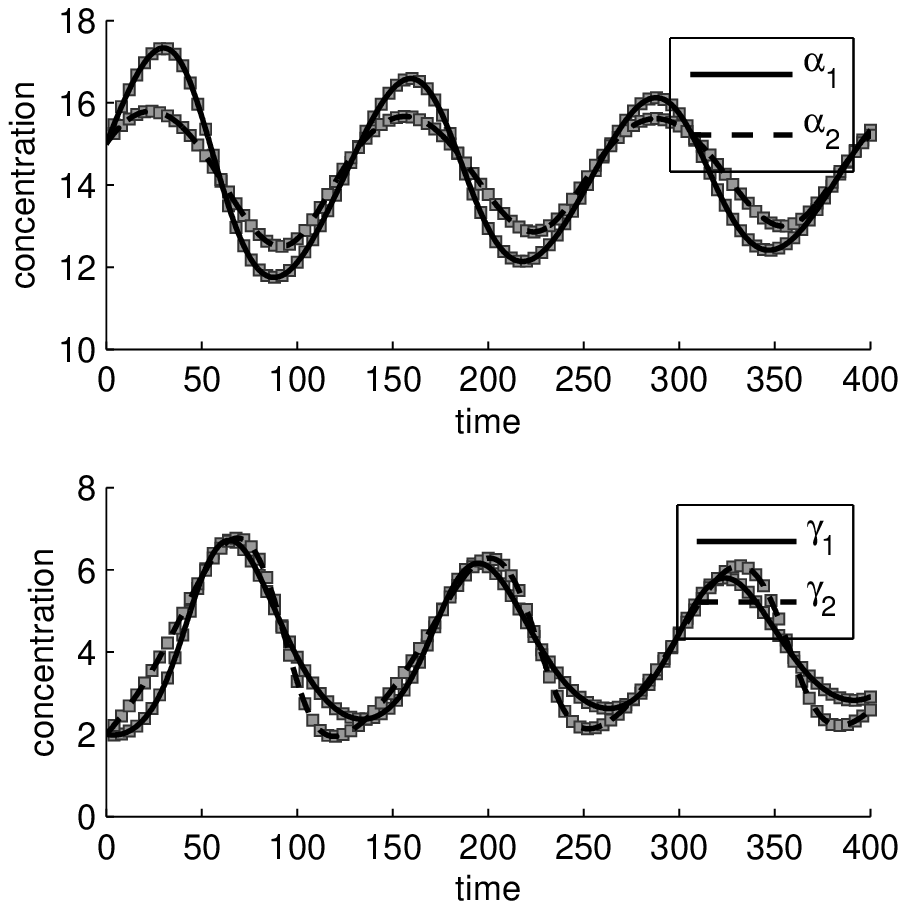}}
\subfigure{\includegraphics[width=0.49\columnwidth]{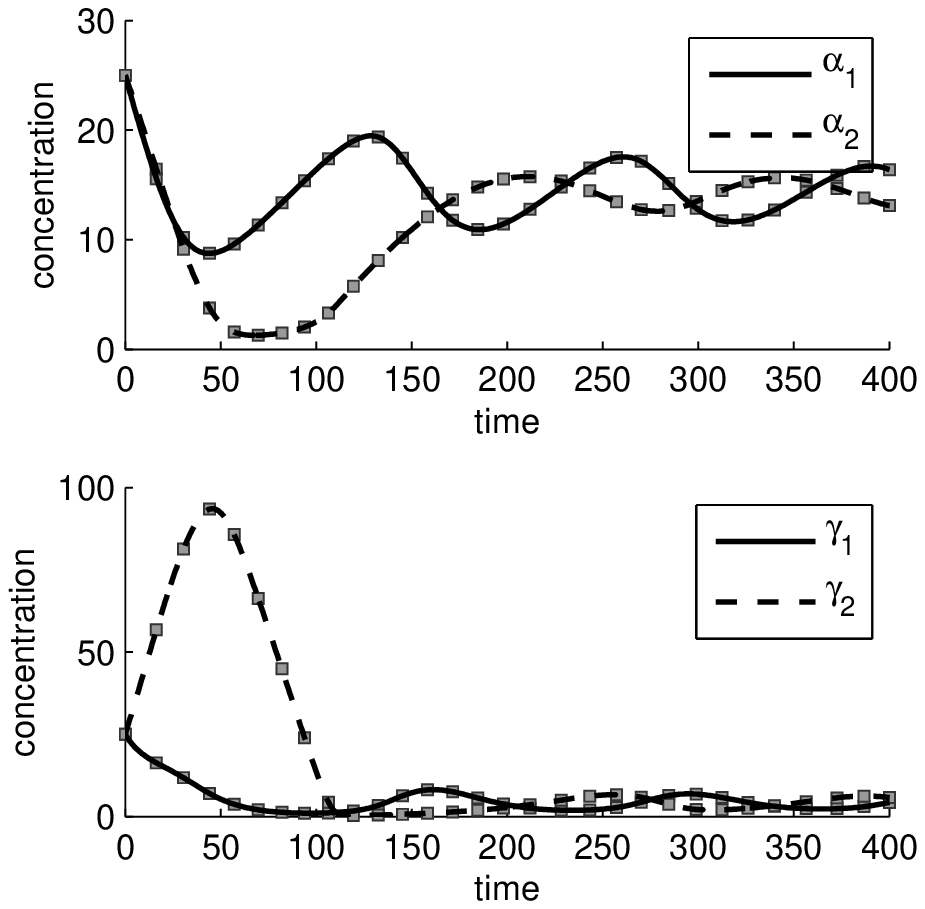}}
\end{center}
\caption{The model functions $\alpha_1,\alpha_2$ and $\gamma_1,\gamma_2$ are shown before the optimization procedure (left) and after the optimization procedure of problem $\mathbf{P}_{\widetilde{\Theta}^1}$ (right) for the glycolytic design setup without the possibility to disturb the system. One square represents one measurement time point.}
\label{glycolyseDESIGNwithoutperturbations1}
\end{figure}

\begin{figure}[ht]
\begin{center}
\subfigure{\includegraphics[width=0.49\columnwidth]{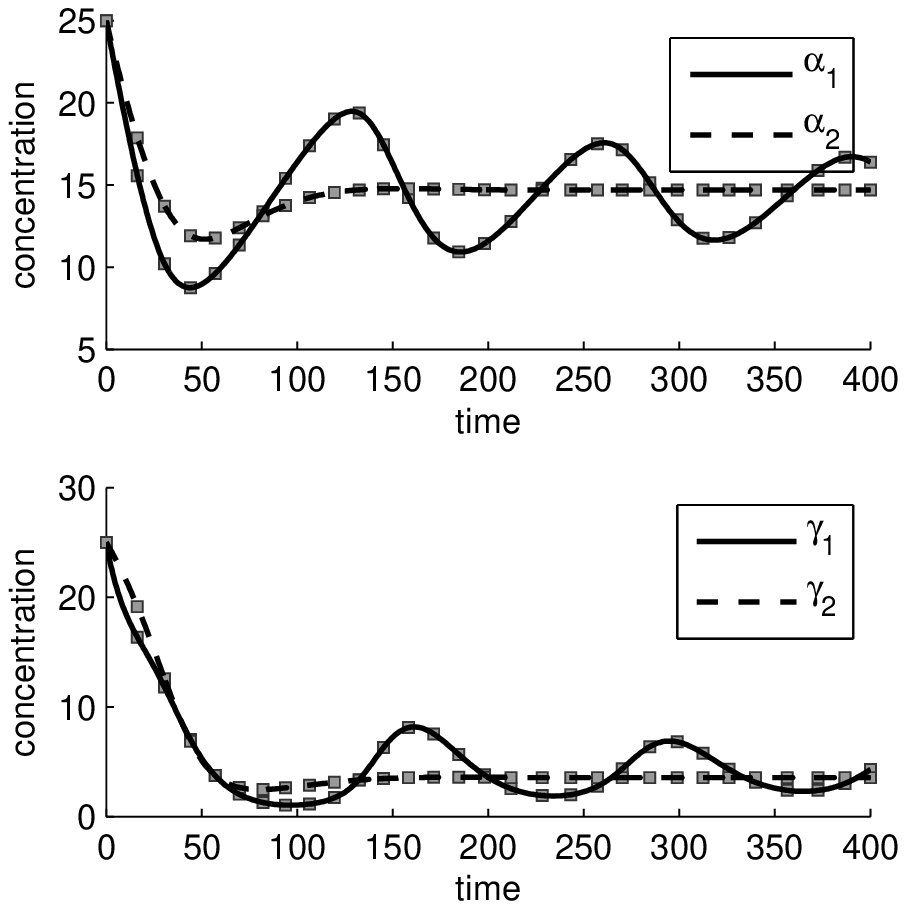}}
\subfigure{\includegraphics[width=0.49\columnwidth]{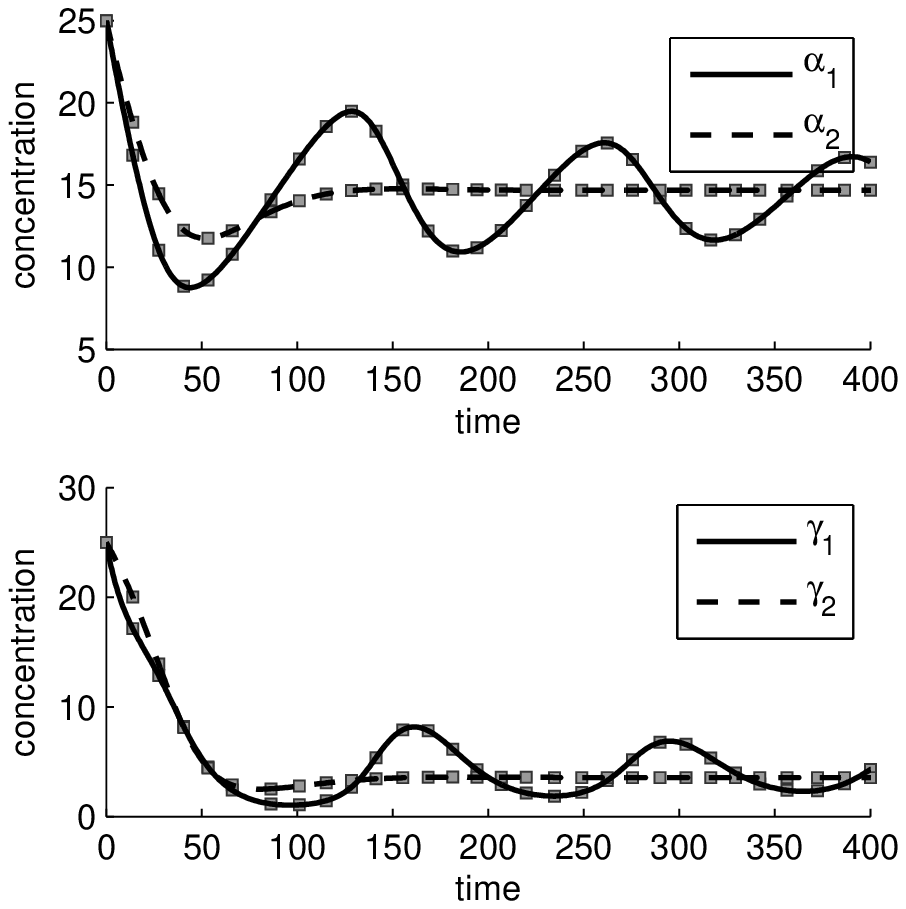}}
\end{center}
\caption{The model functions $\alpha_1,\alpha_2$ and $\gamma_1,\gamma_2$ are shown for the same solution design as for problem $\mathbf{P}_{\widetilde{\Theta}^1}$ \corrFIVEDOM{after the next robustification step} (left) and for the final design (right) for the glycolytic design setup without the possibility to disturb the system. One square represents one measurement time point.}
\label{glycolyseDESIGNwithoutperturbations2}
\end{figure}

\begin{figure}[ht]
\begin{center}
\psfrag{label1}{\parbox{4cm}{\mbox{} \\ \mbox{} \\ $N$}}
\psfrag{label2}{\parbox{4cm}{\vspace{-0.3cm}\hspace{-0.3cm} $\Delta_{RG}$}}
\subfigure{\includegraphics[width=0.45\columnwidth]{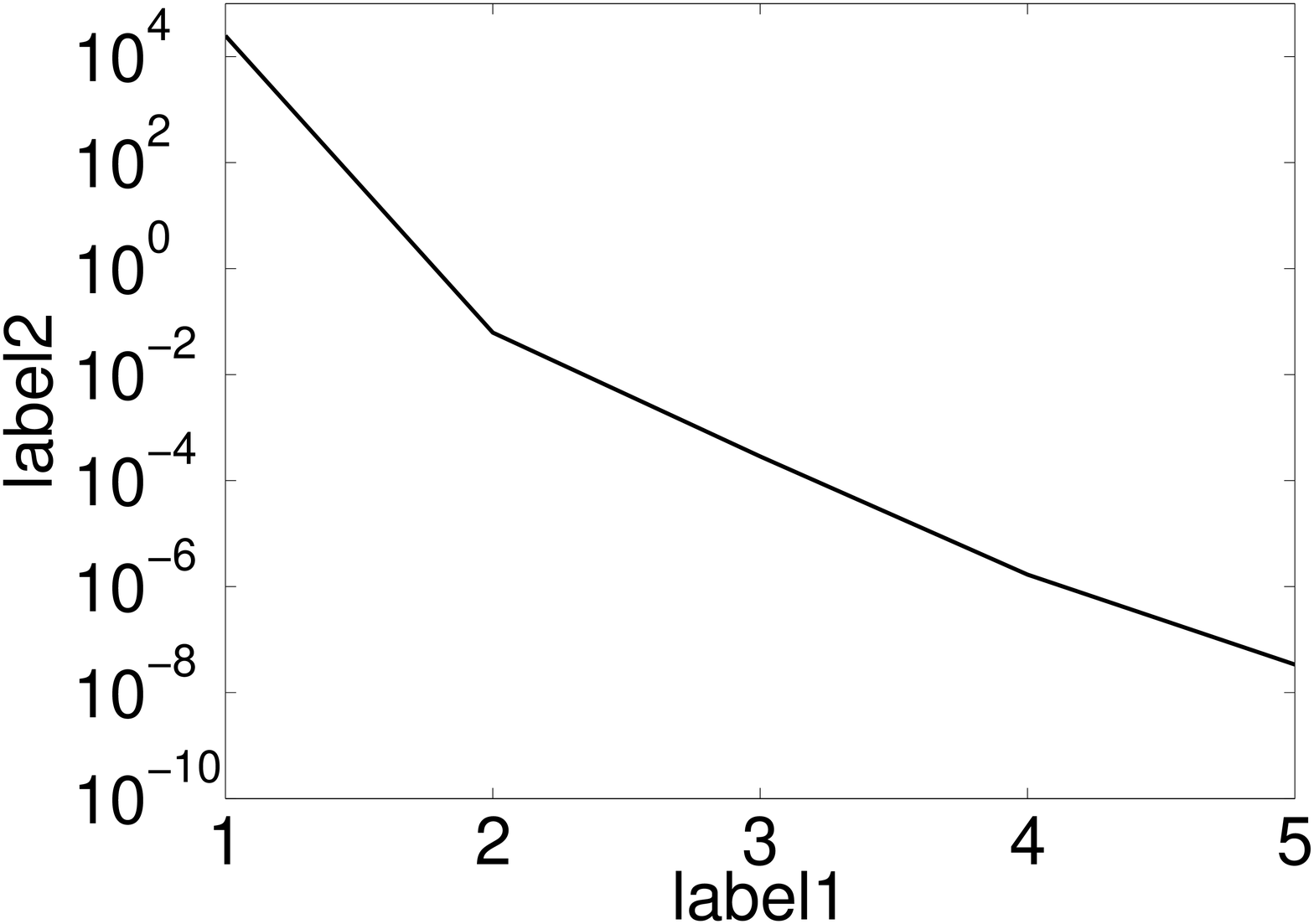}}
\hspace{0.05\columnwidth}
\psfrag{label2}{\parbox{4cm}{\vspace{-0.3cm}\hspace{-1.3cm} objective value of $\mathbf{P}_{\widetilde{\Theta}_N}$}}
\subfigure{\includegraphics[width=0.45\columnwidth]{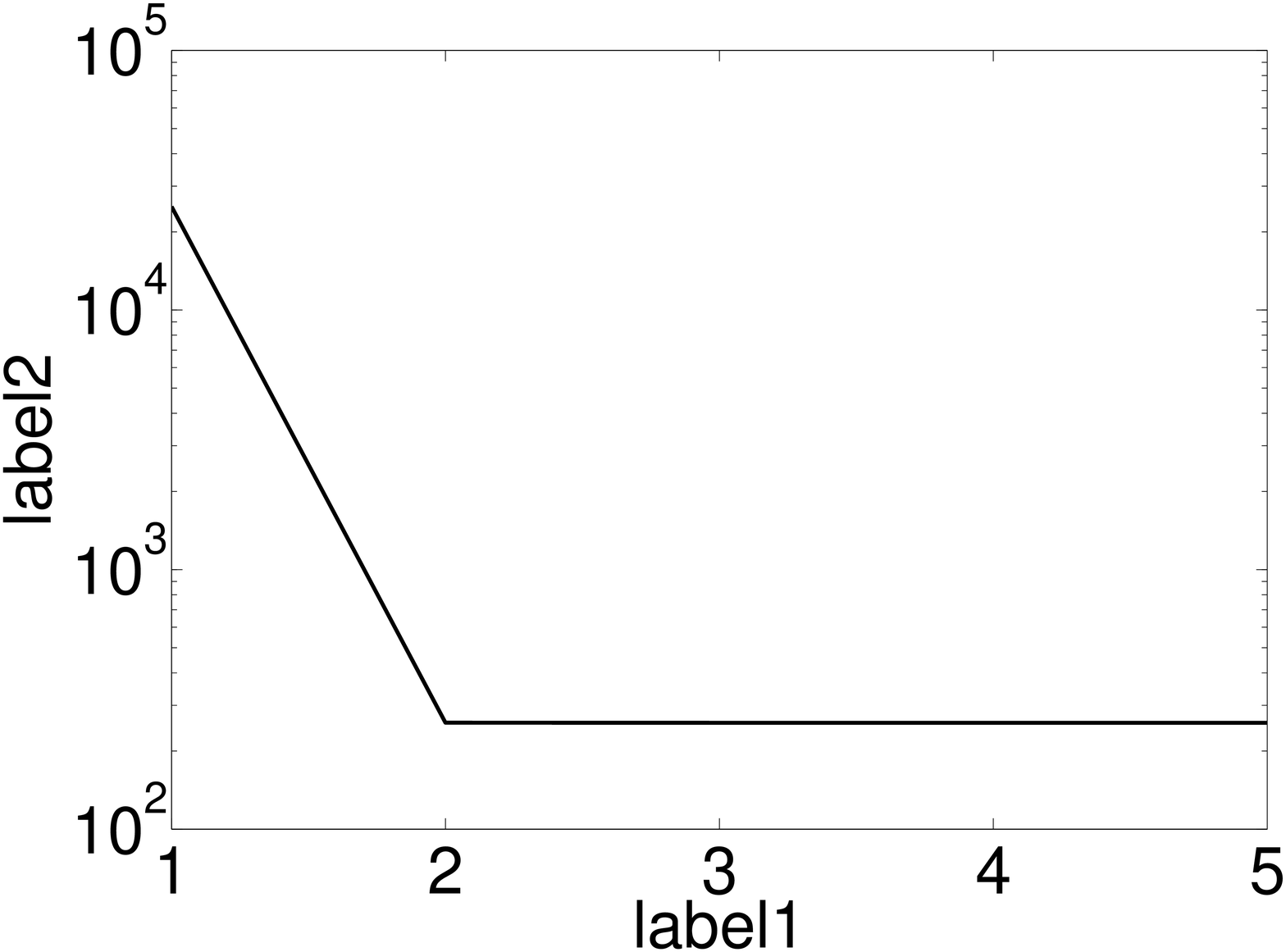}}
\end{center}
\caption{In the left figure the robustification gap $\Delta_{RG}$ is plotted versus the number of iterations $N$ of Algorithm \ref{relaxationAlgorithm} and in the right figure the objective value of problem $\mathbf{P}_{\widetilde{\Theta}^N}$ is shown for the glycolytic design setup without the possibility to disturb the system.}
\label{glycolyseDESIGNwithoutperturbations3}
\end{figure}

\begin{figure}[ht]
\begin{center}
\subfigure{\includegraphics[width=0.49\columnwidth]{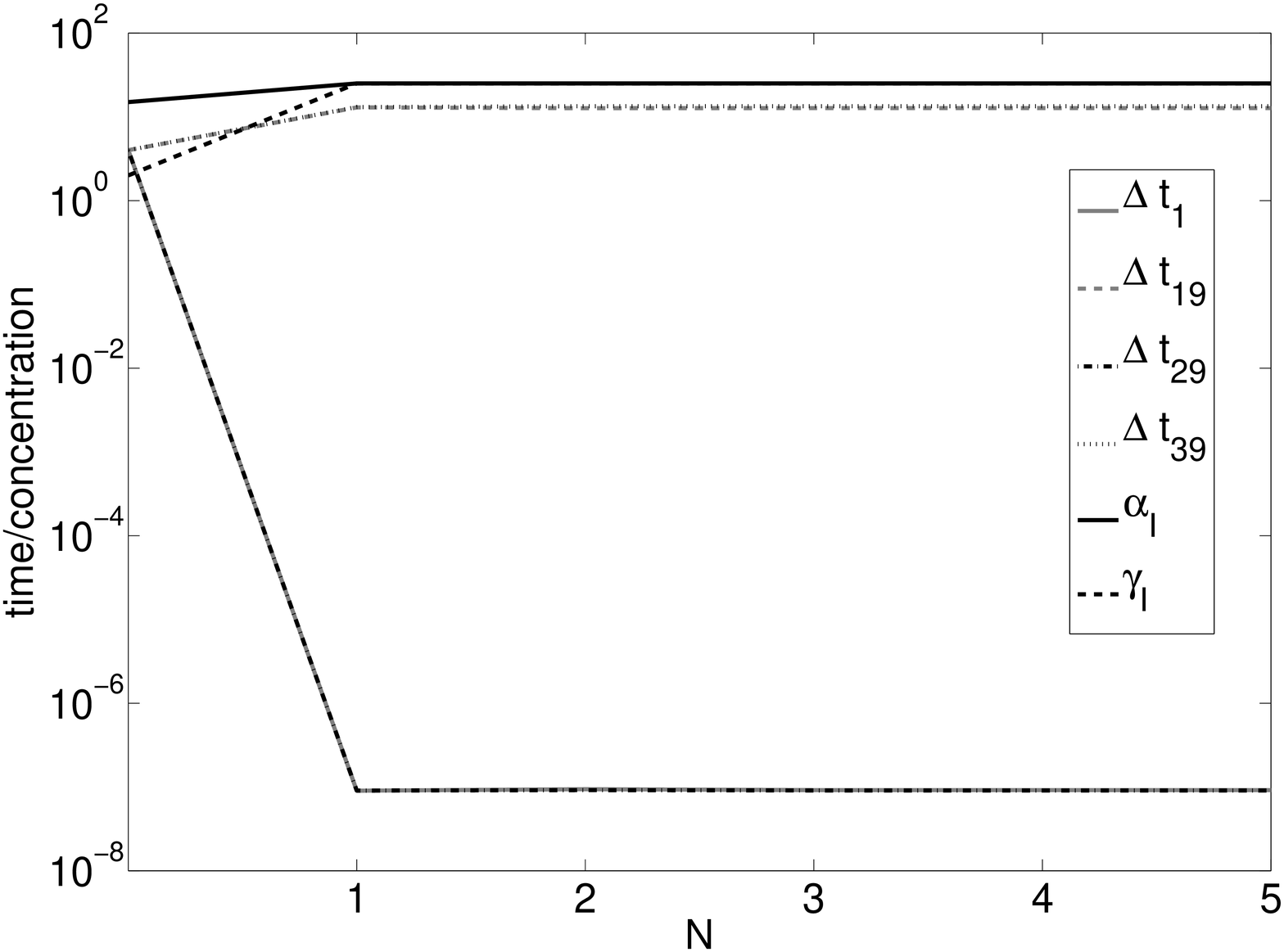}}
\subfigure{\includegraphics[width=0.49\columnwidth]{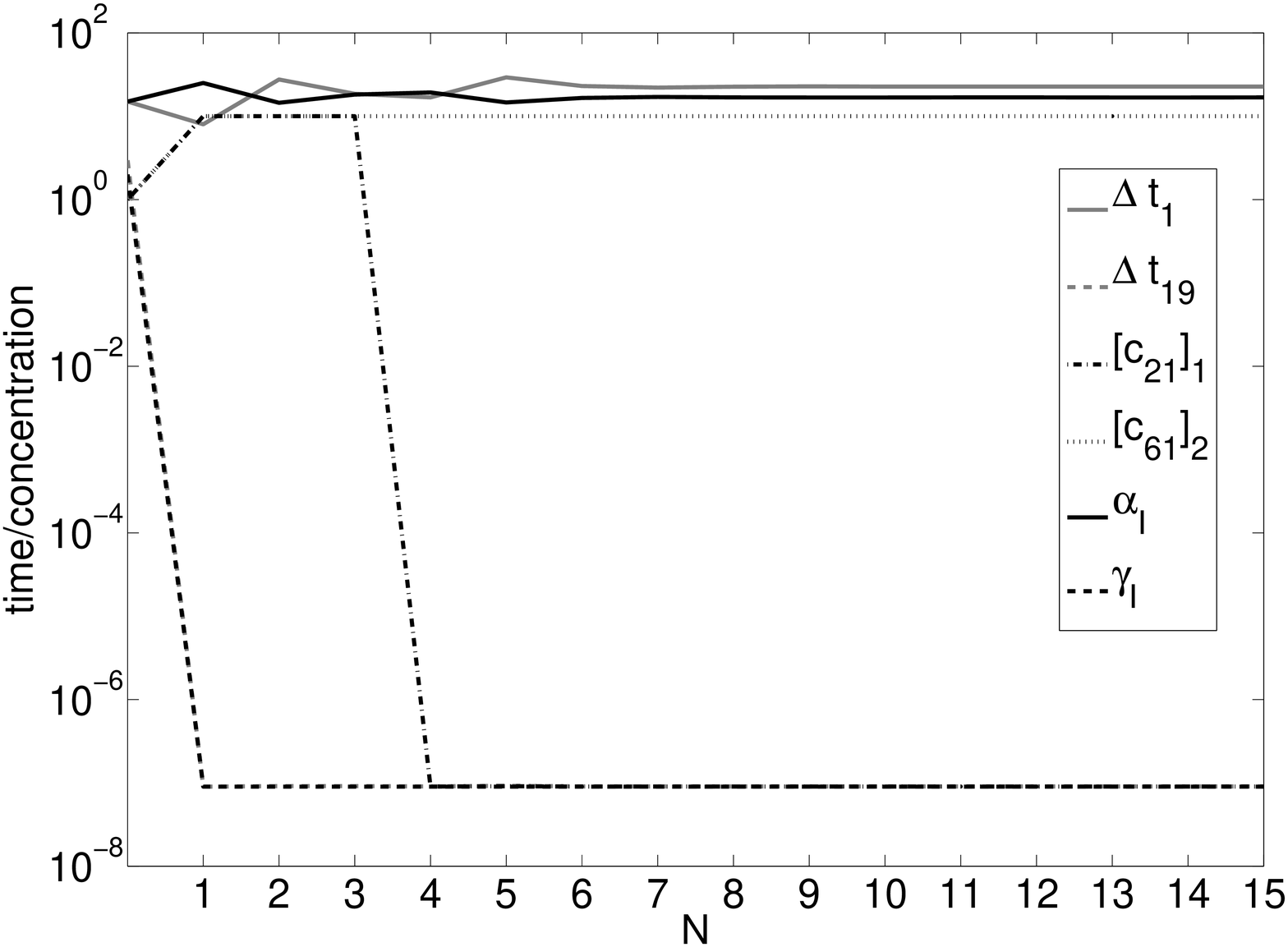}}
\end{center}
\caption{A selection of design variables as solution\corrFIVEDOM{s} of problem $\mathbf{P}_{\widetilde{\Theta}^{N}}$ for the glycolytic design setup without the possibility to disturb the system (left) and with the possibility to disturb the system (rigth) are shown.}
\label{glycolyseDESIGNwithoutperturbations4}
\end{figure}

\begin{figure}[ht]
\begin{center}
\subfigure{\includegraphics[width=0.49\columnwidth]{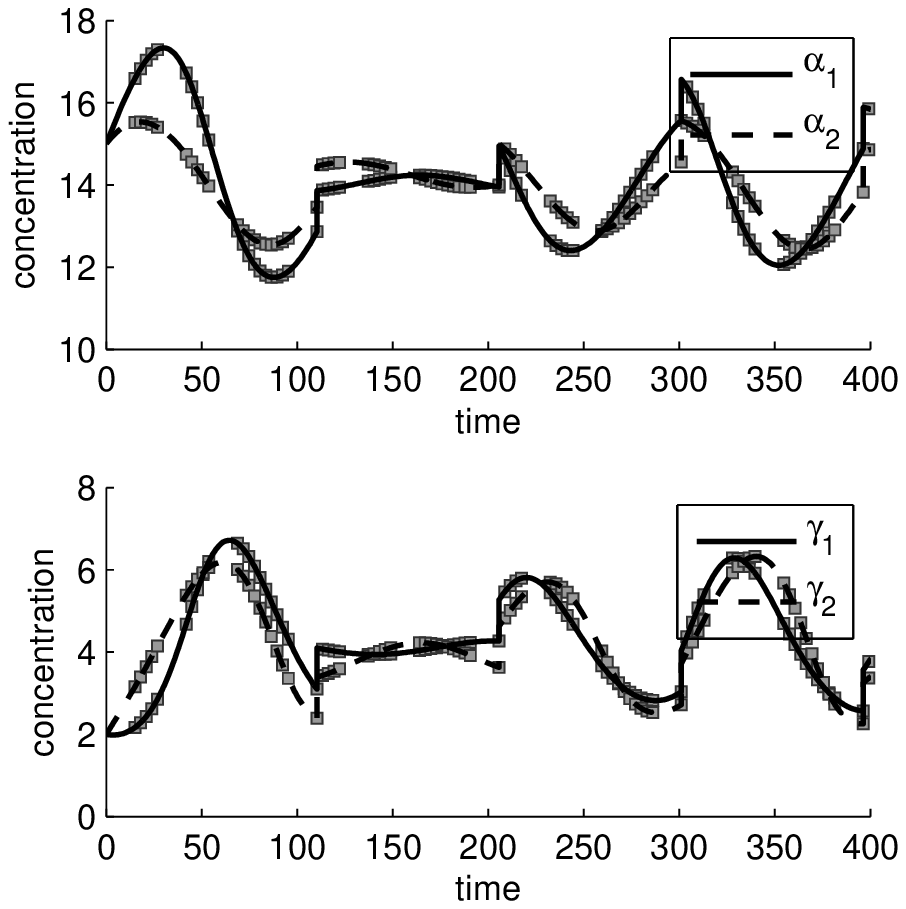}}
\subfigure{\includegraphics[width=0.49\columnwidth]{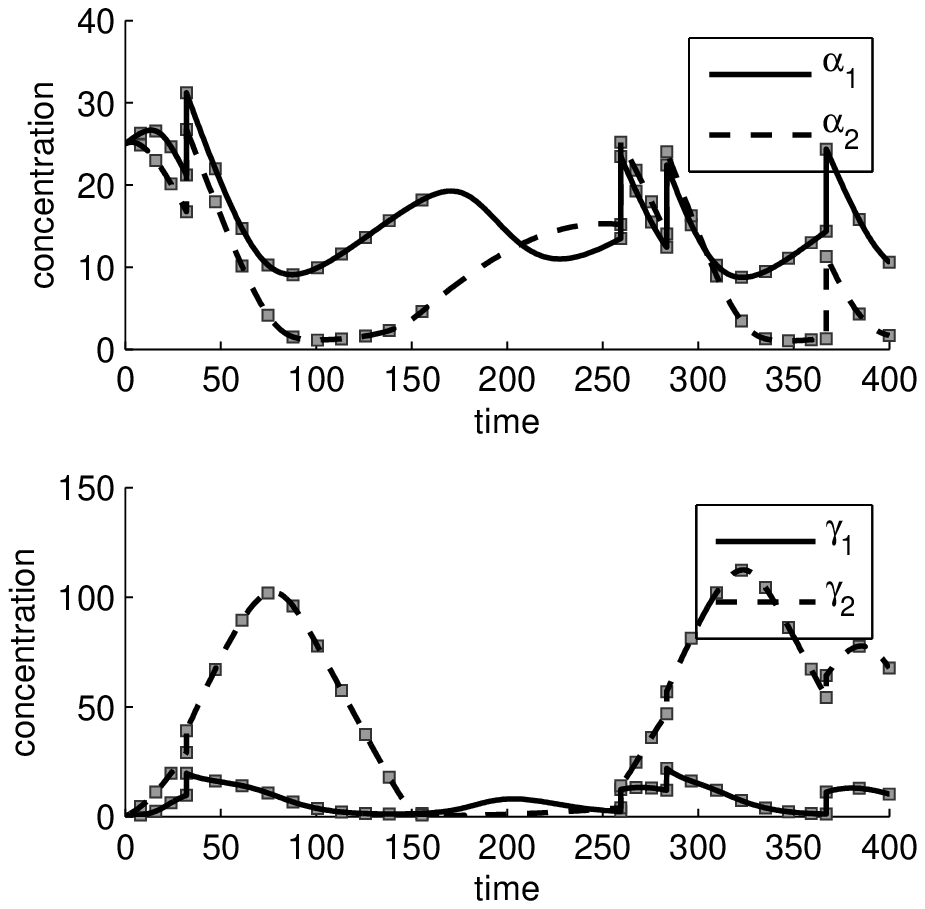}}
\end{center}
\caption{The model functions $\alpha_1,\alpha_2$ and $\gamma_1,\gamma_2$ are shown before the optimization procedure (left) and after the optimization procedure of problem $\mathbf{P}_{\widetilde{\Theta}^1}$ (right) for the glycolytic design setup with the possibility to disturb the system. One square represents one measurement time point.}
\label{glycolyseDESIGNwithperturbations1}
\end{figure}

\begin{figure}[ht]
\begin{center}
\subfigure{\includegraphics[width=0.49\columnwidth]{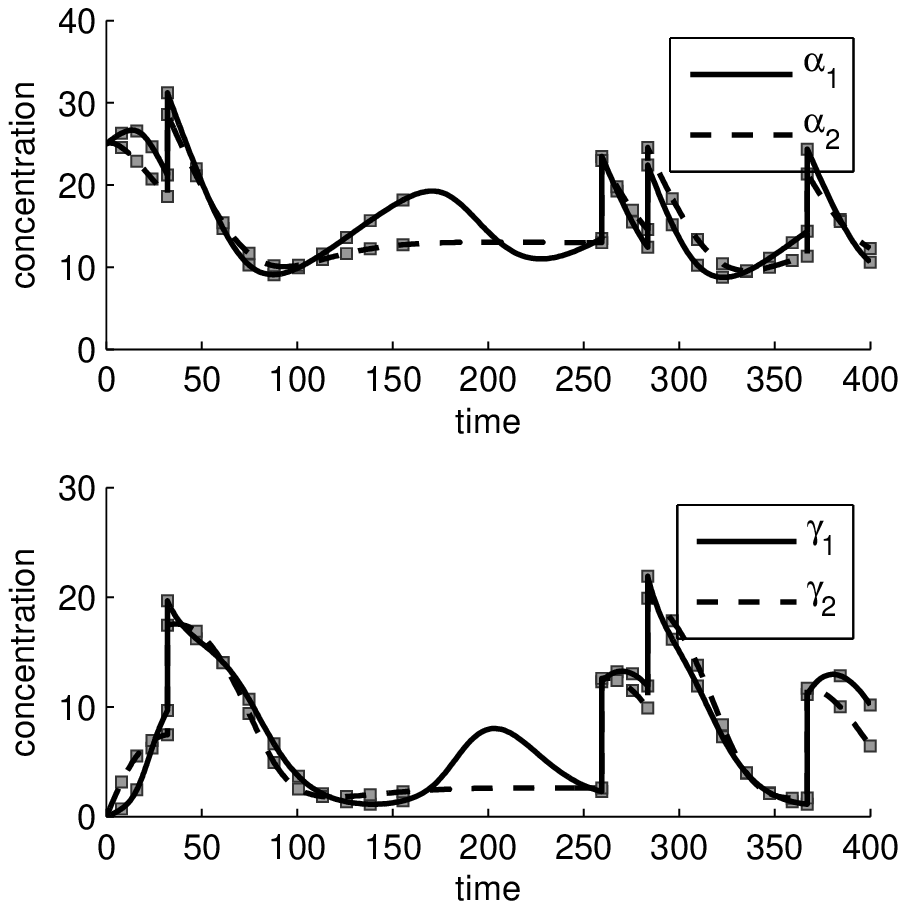}}
\subfigure{\includegraphics[width=0.49\columnwidth]{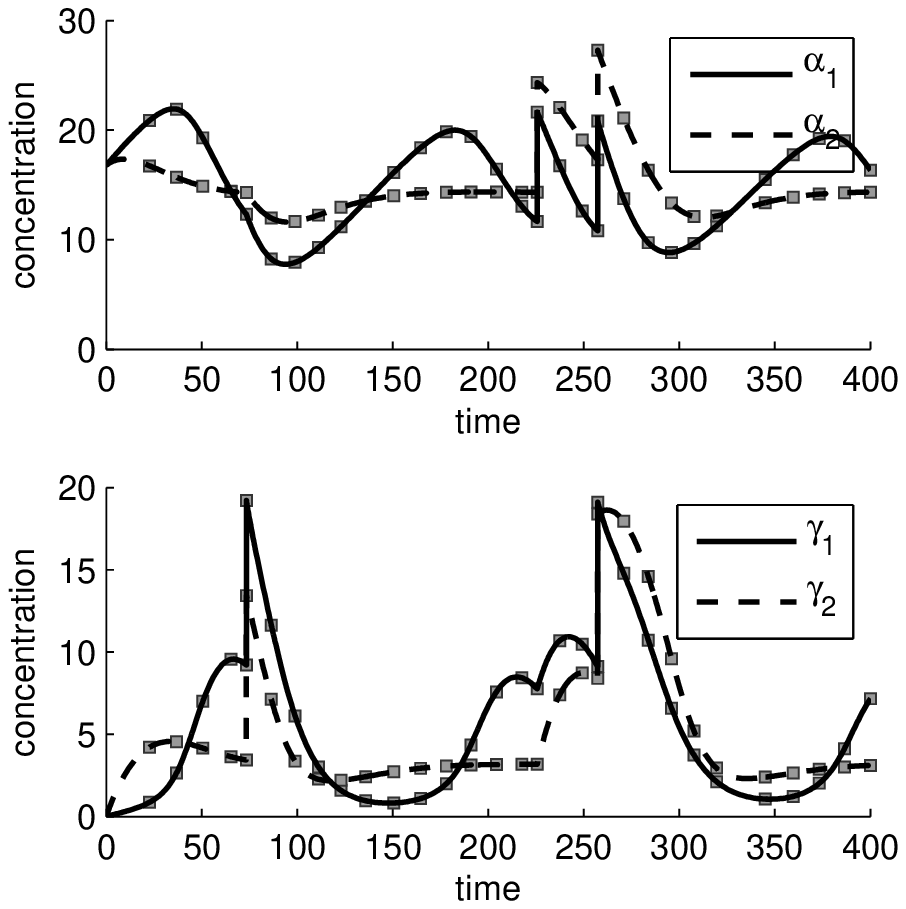}}
\end{center}
\caption{The model functions $\alpha_1,\alpha_2$ and $\gamma_1,\gamma_2$ are shown for the same solution design as for problem $\mathbf{P}_{\widetilde{\Theta}^1}$ \corrFIVEDOM{after the next robustification step} (left) and for the final design (right) for the glycolytic design setup with the possibility to disturb the system. One square represents one measurement time point.}
\label{glycolyseDESIGNwithperturbations2}
\end{figure}

\begin{figure}[ht]
\begin{center}
\psfrag{label1}{\parbox{4cm}{\mbox{} \\ \mbox{} \\ $N$}}
\psfrag{label2}{\parbox{4cm}{\vspace{-0.3cm}\hspace{-0.3cm} $\Delta_{RG}$}}
\subfigure{\includegraphics[width=0.45\columnwidth]{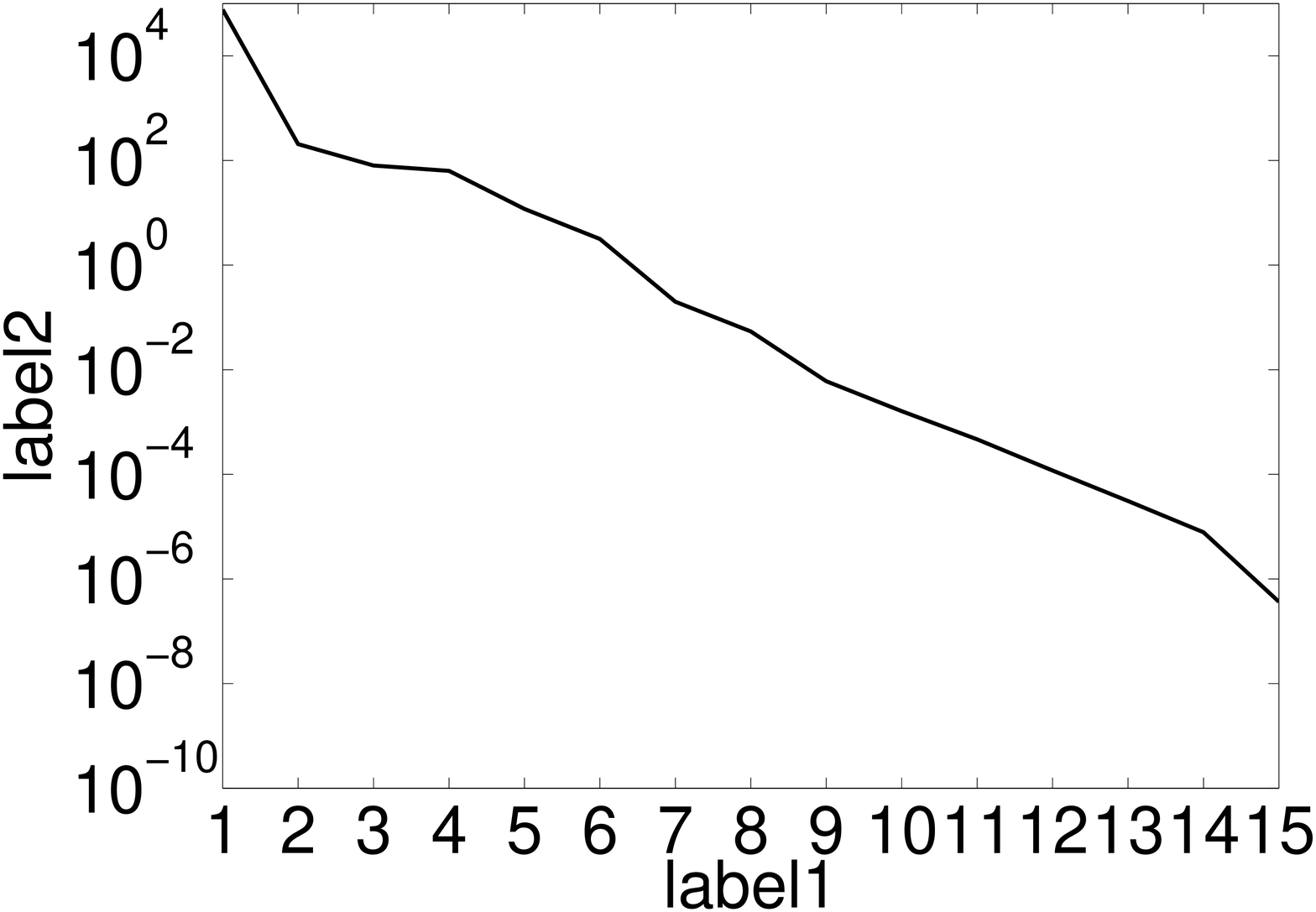}}
\hspace{0.05\columnwidth}
\psfrag{label2}{\parbox{4cm}{\vspace{-0.3cm}\hspace{-1.3cm} objective value of $\mathbf{P}_{\widetilde{\Theta}_N}$}}
\subfigure{\includegraphics[width=0.45\columnwidth]{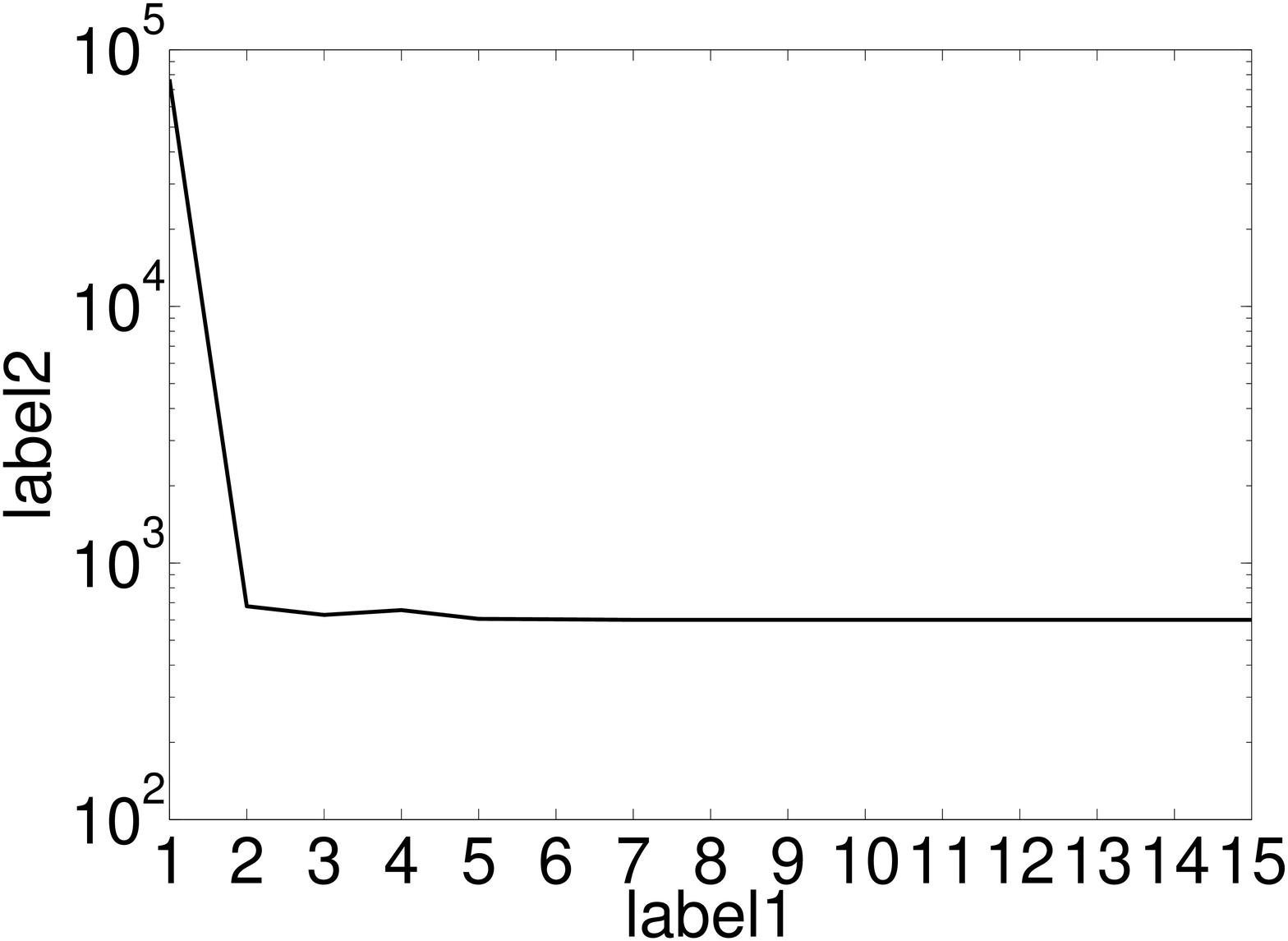}}
\end{center}
\caption{In the left figure the robustification gap $\Delta_{RG}$ is plotted versus the number of iterations $N$ of Algorithm \ref{relaxationAlgorithm} and in the right figure the objective value of problem $\mathbf{P}_{\widetilde{\Theta}^N}$ is shown for the glycolytic design setup with the possibility to disturb the system.}
\label{glycolyseDESIGNwithperturbations3}
\end{figure}

In the first test case for model discrimination we implemented the following models for glycolytic oscillations as described in \cite{Goldbeter1996}.\\\\
Model 1 is an allosteric enzyme model with positive feedback under cooperativity and linear product sink. The differential equations for model 1 are given by
\begin{equation*}
\begin{split}
  \dfrac{d\alpha_1}{dt}&=\nu -\sigma \phi(\alpha_1,\gamma_1 ), \\
 \dfrac{d\gamma_1}{dt}&=q_1 \sigma \phi(\alpha_1,\gamma_1 ) -k_s \gamma_1, \\
  \phi(\alpha_1,\gamma_1 )&=\dfrac{\alpha_1(1+\alpha_1)(1+\gamma_1)^2}{L_1+(1+\alpha_1)^2(1+\gamma_1)^2}. 
\end{split}
\end{equation*}\mbox{}\\\\
Model 2 is an allosteric model with positive feedback in the absence of cooperativity and the product sink is represented by Michaelis-Menten kinetics. The differential equations for Model 2 are given by
\begin{equation*}
\begin{split}
  \dfrac{d\alpha_2}{dt}&=\nu - \phi(\alpha_2,\gamma_2 ), \\
 \dfrac{d\gamma_2}{dt}&=q_2 \phi(\alpha_2,\gamma_2) -\dfrac{r_s \gamma_2}{\mu+\gamma_2}, \\
  \phi(\alpha_2,\gamma_2 )&=\dfrac{\alpha_2 (1+\gamma_2)}{L_2+(1+\alpha_2)(1+\gamma_2)}. 
\end{split}
\end{equation*}\mbox{}\\\\
$\alpha_{1,2}$ denotes the species concentration of the substrate $\gamma_{1,2}$ that of the product.\\ 
For both models the inflow parameter $\nu$ is the same and fixed to the value $\nu=0.22$. It represents the inflow of substrate to the experimental system, a CSTR (continuously stirred tank reactor).\\ The parameters $\sigma$, $q_1$, $k_s$ and $L_1$ of model 1 are regarded as known. Their values are given in Table \ref{paraValueNumExampl1}, the parameters $q_2$, $r_s$, $\mu$ and $L_2$ of model 2 are regarded as unknown and subject to robustification. For the permitted parameter range see Table \ref{paraValueNumExampl1}. 
\begin{table}[h]
\begin{center}
\begin{tabular}{||c|c|c|c||c|c|c|c||}
\hline 
\multicolumn{4}{||c||}{Model 1}  & \multicolumn{4}{c||}{Model 2} \tabularnewline
\hline
 $\sigma$ & $q_1$ & $k_s$ & $L_1$ & $q_2$ & $r_s$ & $\mu$ & $L_2$\tabularnewline
 0.92 & 2.01 & 0.11 & 17206.10 & $[10^{-7},100]$ & $[10^{-7},100]$ & $[10^{-7},100]$ & $[100,300]$\tabularnewline
\hline
\end{tabular}
\end{center}
\caption{Parameter values for the glycolytic oscillation models.}
\label{paraValueNumExampl1}
\end{table}\\ 
For simplicity we consider the homoscedastic case with equal variances, i.e. $v_1 = v_2 =\sigma^2$. In this case $\mathcal{I}(P_1:P_2,\mathcal{O}_1)$ reduces to,
\begin{equation}
 \mathcal{I}(P_1:P_2,\mathcal{O}_1)=\sum_{i=1}^n  \corrFIVEDOM{\mathcal{H}'} (t_i)\corrFIVEDOM{\widetilde{\mathcal{H}'}}(c_i)\left((\alpha^i_1-\alpha^i_2)^2+(\gamma^i_1-\gamma^i_2)^2\right).
\end{equation}
For this test case the homotopy strategy as presented in Section \ref{HomotopyStrategy} is only applied if the robustification gap $\Delta_{RG}<0.1$, then the successive problem $\mathbf{P}_{\widetilde{\Theta}^{N+1}}$ is calculated by \corrFIVEDOM{use of} the homotopy strategy with $30$ homotopy steps, i.e. $\kappa_h=h/30$, $h\in\{1,...,30\}$. Otherwise the problem $\mathbf{P}_{\widetilde{\Theta}^{N+1}}$ is \corrFIVEDOM{solved} without homotopy strategy.
For each subsequent problem $\mathbf{P}_{\widetilde{\Theta}^{N+1}}$ the solution of problem $\mathbf{P}_{\widetilde{\Theta}^{N}}$ is used as initial guess.\\\\
We first present a robust design without the possibility to disturb the system by adding species at later time points.\\ 
The design is calculated within a fixed time window i.e. $T_{\mathrm{end}}=400$. $100$ equally spaced possible measurement points are defined in the initial state of the optimization procedure, the distance vector $\Delta t$ between the time points is subject to design and each entry is restricted to $\Delta t_i\in[10^{-7},10^{19}],$ $i\in \{1,...,100\}$. The disturbance vectors $c_i$ are set to $c_i=0$, $i\in\{1,...,99\}$ and are fixed to model the fact that no species disturbance is allowed.\\ 
The initial species concentrations which are also subject to experimental design are restricted 
to $\alpha_I \in [10^{-7},25]$ and $\gamma_I \in [10^{-7},25]$. The initial values were set to $\alpha_I=15$ and $\gamma_I=2$. The parameters of the switching functions $\corrFIVEDOM{\mathcal{H}'}(t_i)$ are chosen as $a_1=20.0$ and $b_1=10.0$. The parameters of the switching functions $\corrFIVEDOM{\widetilde{\mathcal{H}'}}(c_i)$ are chosen as $a_2=0.05$ and $b_2=0.025$. The algorithmic settings are summarized in Table \ref{numericSettingsExample1}.
\begin{table}[h]
\begin{center}
\begin{tabular}{||c|c|c||c|c||}
\hline 
\multicolumn{3}{||c||}{Optimization settings}  & \multicolumn{2}{c||}{Integrator settings} \tabularnewline
\hline
 $P$ & $\delta$ &  IPOPT-tol: Step 1./Step 2. & relTol/absTol & relTolSens/absTolSens \tabularnewline
 $5$ & $10^{-6}$ & $10^{-10}$/$10^{-8}$ & $10^{-12}$/$10^{-12}$ &  $10^{-12}$/$10^{-12}$ \tabularnewline
\hline
\end{tabular}
\end{center}
\caption{On the left \corrFIVEDOM{hand side} the optimization settings are listed comprising the IPOPT stopping tolerances for Step 1. and Step 2. of Algorithm \ref{relaxationAlgorithm} and on the right \corrFIVEDOM{hand side} the integration tolerances for the nominal trajectory and the first \corrFIVEDOM{order} sensitivities. We use the IPOPT option  ``honor\_original\_bounds=no'' for Step 1. and Step 2. of Algorithm \ref{relaxationAlgorithm}.}
\label{numericSettingsExample1}
\end{table}\\
A plot of the functions $\alpha_1,\alpha_2$ and $\gamma_1,\gamma_2$ \corrFIVEDOM{in} the initial state and for the solution of problem $\mathbf{P}_{\widetilde{\Theta}^1}$ are shown in Figure \ref{glycolyseDESIGNwithoutperturbations1}. A plot for the same functions with the same solution design as for problem $\mathbf{P}_{\widetilde{\Theta}^1}$ \corrFIVEDOM{after the next robustification step} is shown in Figure \ref{glycolyseDESIGNwithoutperturbations2}. The final design is also shown in Figure \ref{glycolyseDESIGNwithoutperturbations2}.
A plot of the robustification gap $\Delta_{RG}$ and as well for the objective value of problem $\mathbf{P}_{\widetilde{\Theta}^{N}}$ for each iteration $N$ of Algorithm \ref{relaxationAlgorithm} are shown in Figure \ref{glycolyseDESIGNwithoutperturbations3}. A selection of design variables as solution\corrFIVEDOM{s} of problem $\mathbf{P}_{\widetilde{\Theta}^{N}}$ is shown in Figure \ref{glycolyseDESIGNwithoutperturbations4}(left). 
\\\\
In a second scenario we additionally allow for species perturbations. \corrFIVEDOM{In this} new scenario at the $21$th $41$th, $61$th  and $81$th measurement time points, the system can get disturbed by additional species quantities. The free vectors $c_i$, $i\in\{21,41,61,81\}$ are \corrFIVEDOM{constrained by} $c_i\in[10^{-7},10]$. The initial values are set to $c_i=1$. The remaining conditions are as before, however we change the time vector bound constraints for $i\in\{1,6,11,21,26,31,41,46,51,61,66,71,81\}$ to $\Delta t_i\in[8,10^{19}]$ and the initial state to $\Delta t_i=15$. The bounds for the remaining entries are as before, and the remaining measurement time points were equally spaced.\\ 
A plot of the functions $\alpha_1,\alpha_2$ and $\gamma_1,\gamma_2$ \corrFIVEDOM{in} the initial state and for the solution of problem $\mathbf{P}_{\widetilde{\Theta}^1}$ are shown in Figure \ref{glycolyseDESIGNwithperturbations1}. A plot for the same functions with the same solution design as for problem $\mathbf{P}_{\widetilde{\Theta}^1}$ \corrFIVEDOM{after the next robustification step} is shown in Figure \ref{glycolyseDESIGNwithperturbations2}. The final design is also shown in Figure \ref{glycolyseDESIGNwithperturbations2}.
A plot of the robustification gap $\Delta_{RG}$ \corrFIVEDOM{and} the objective value of problem $\mathbf{P}_{\widetilde{\Theta}^{N}}$ for each iteration $N$ of Algorithm \ref{relaxationAlgorithm} are shown in Figure \ref{glycolyseDESIGNwithperturbations3}. A selection of design variables as solution\corrFIVEDOM{s} of problem $\mathbf{P}_{\widetilde{\Theta}^{N}}$ is shown in Figure \ref{glycolyseDESIGNwithoutperturbations4}(right).
\subsection{Discriminating design for two models describing signal sensing in dictyostelium discoideum}
\begin{figure}[ht]
\begin{center}
\subfigure[model 1]{\includegraphics[width=0.2\columnwidth]{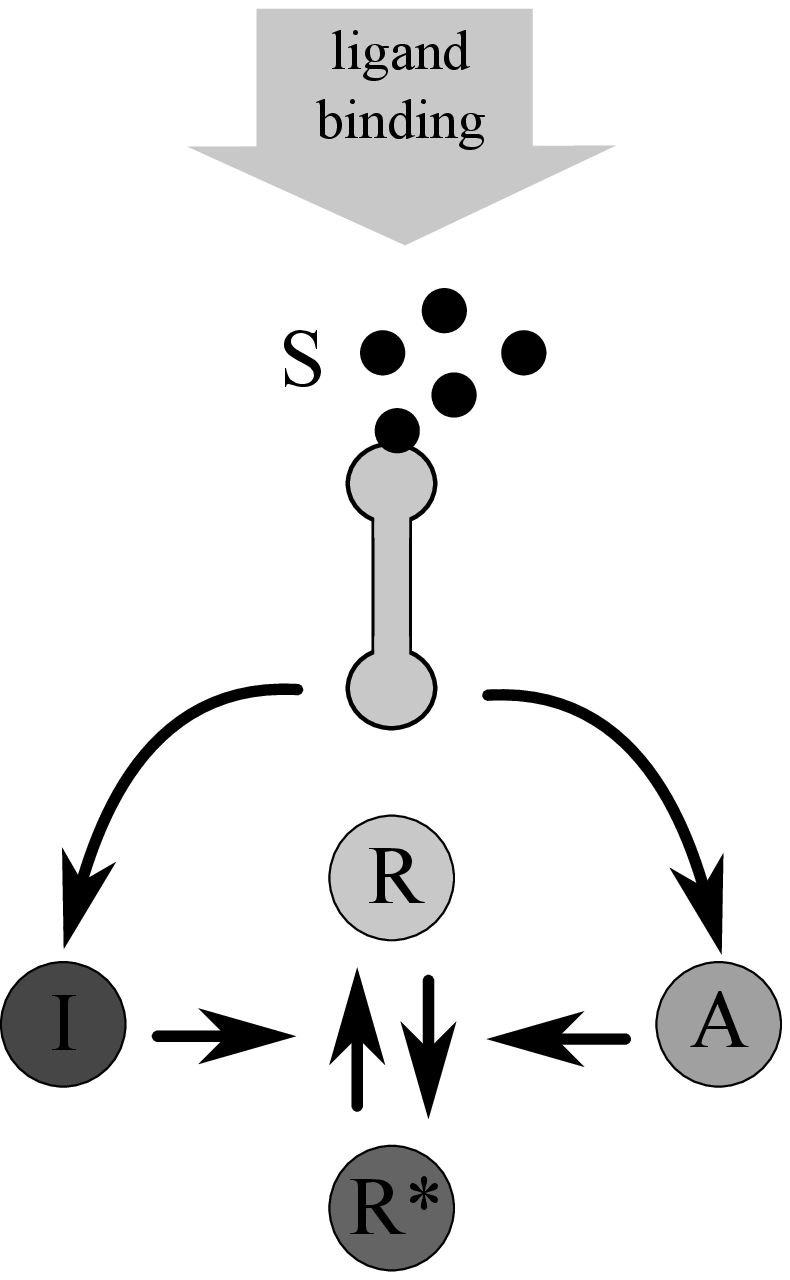}}
\hspace{0.25\columnwidth} 
\subfigure[model 2]{\includegraphics[width=0.2\columnwidth]{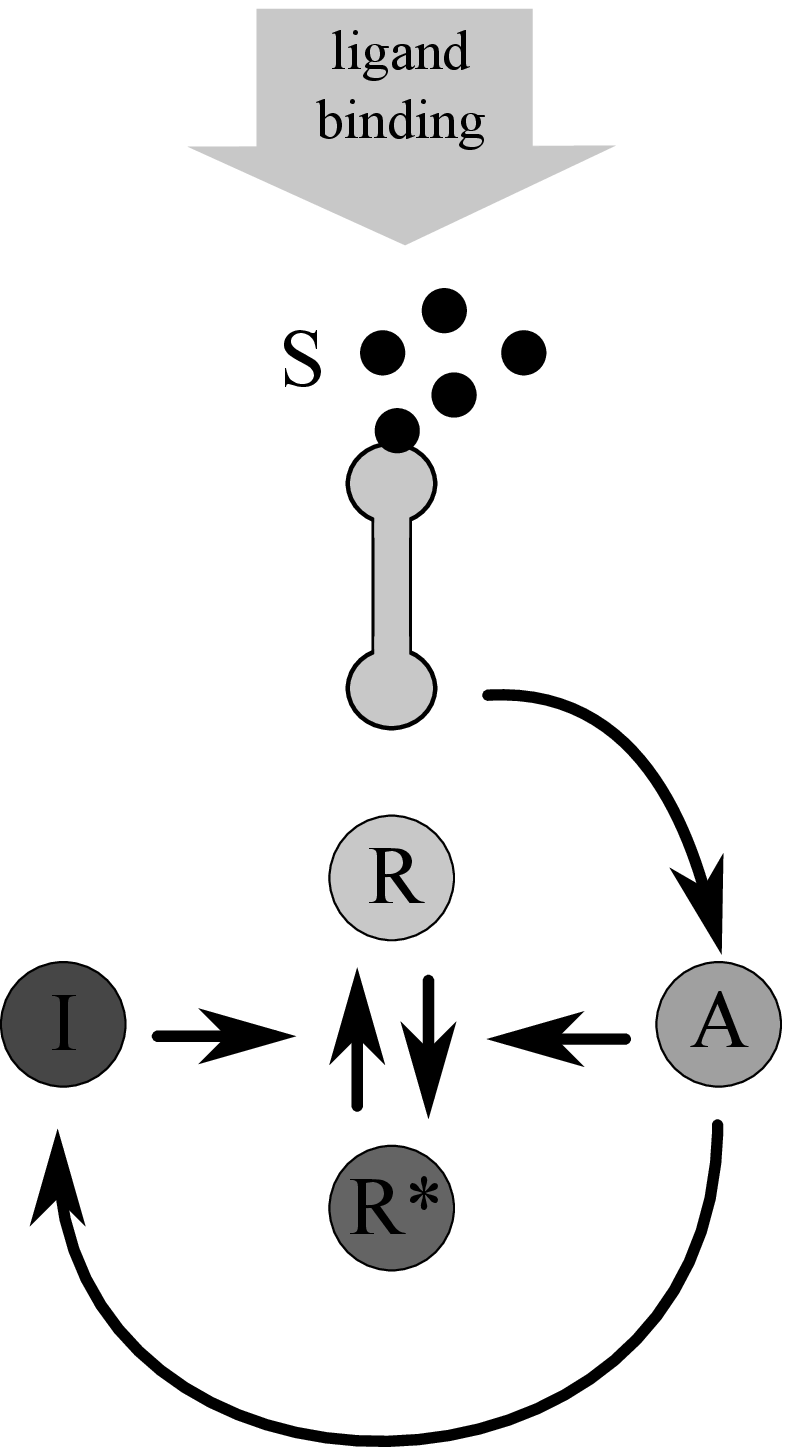}}
\end{center}
\caption{Two models of the signal system of the Dictyostelium amoeba.}
\label{amoebaModel}
\end{figure}

\begin{figure}[ht]
\begin{center}
\subfigure{\includegraphics[width=0.49\columnwidth]{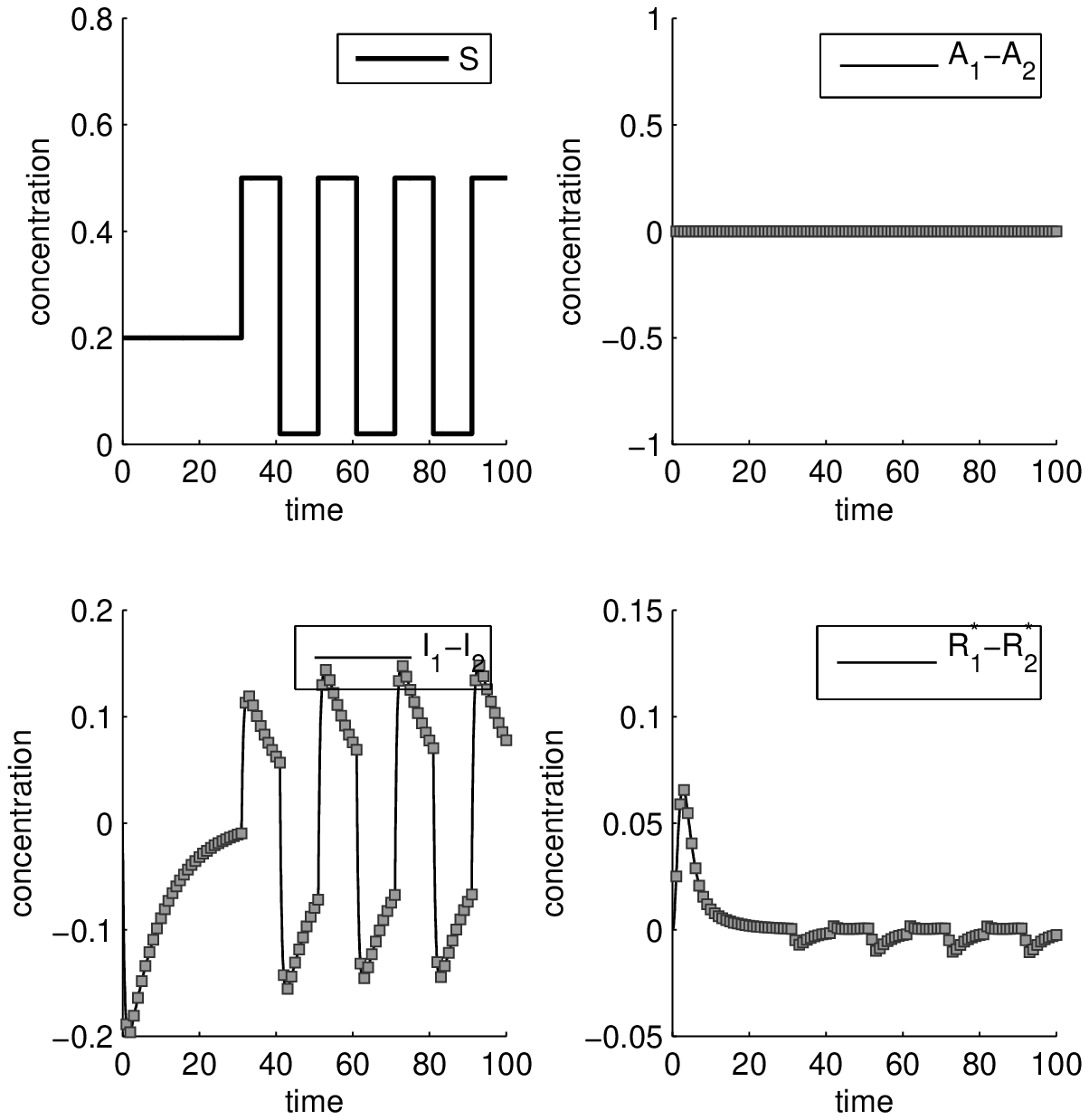}}
\subfigure{\includegraphics[width=0.49\columnwidth]{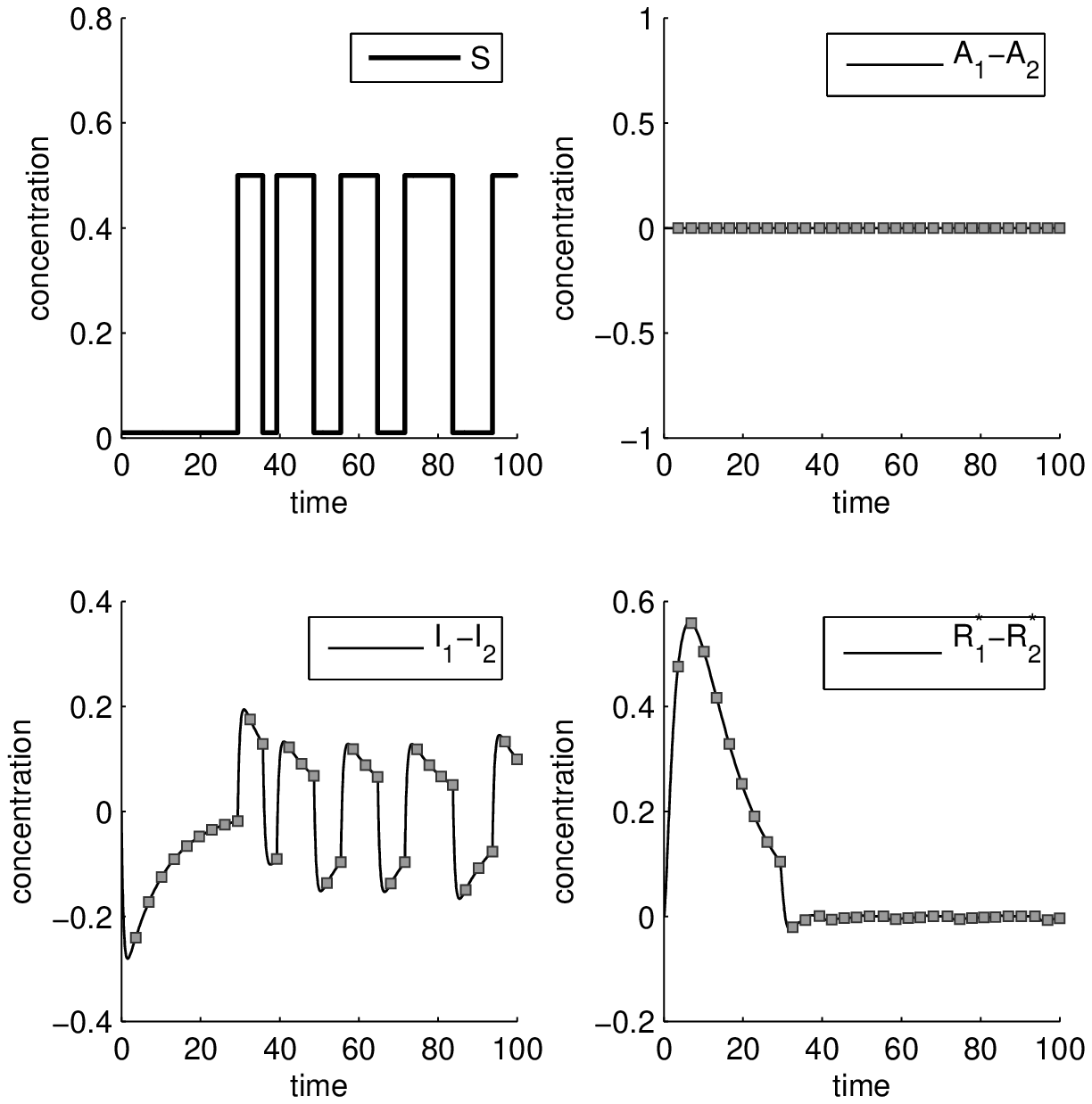}}
\end{center}
\caption{The model \corrFIVEDOM{variable} distance functions $(S_1-S_2)$, $(A_1-A_2)$, $(I_1-I_2)$ and $(R_1-R_2)$ are shown before the optimization procedure (left) and after the optimization procedure of problem $\mathbf{P}_{\widetilde{\Theta}^1}$ (right) for two models describing signal sensing in dictyostelium discoideum. One square represents one measurement time point.}
\label{amoebaDESIGN1}
\end{figure}

\begin{figure}[ht]
\begin{center}
\subfigure{\includegraphics[width=0.49\columnwidth]{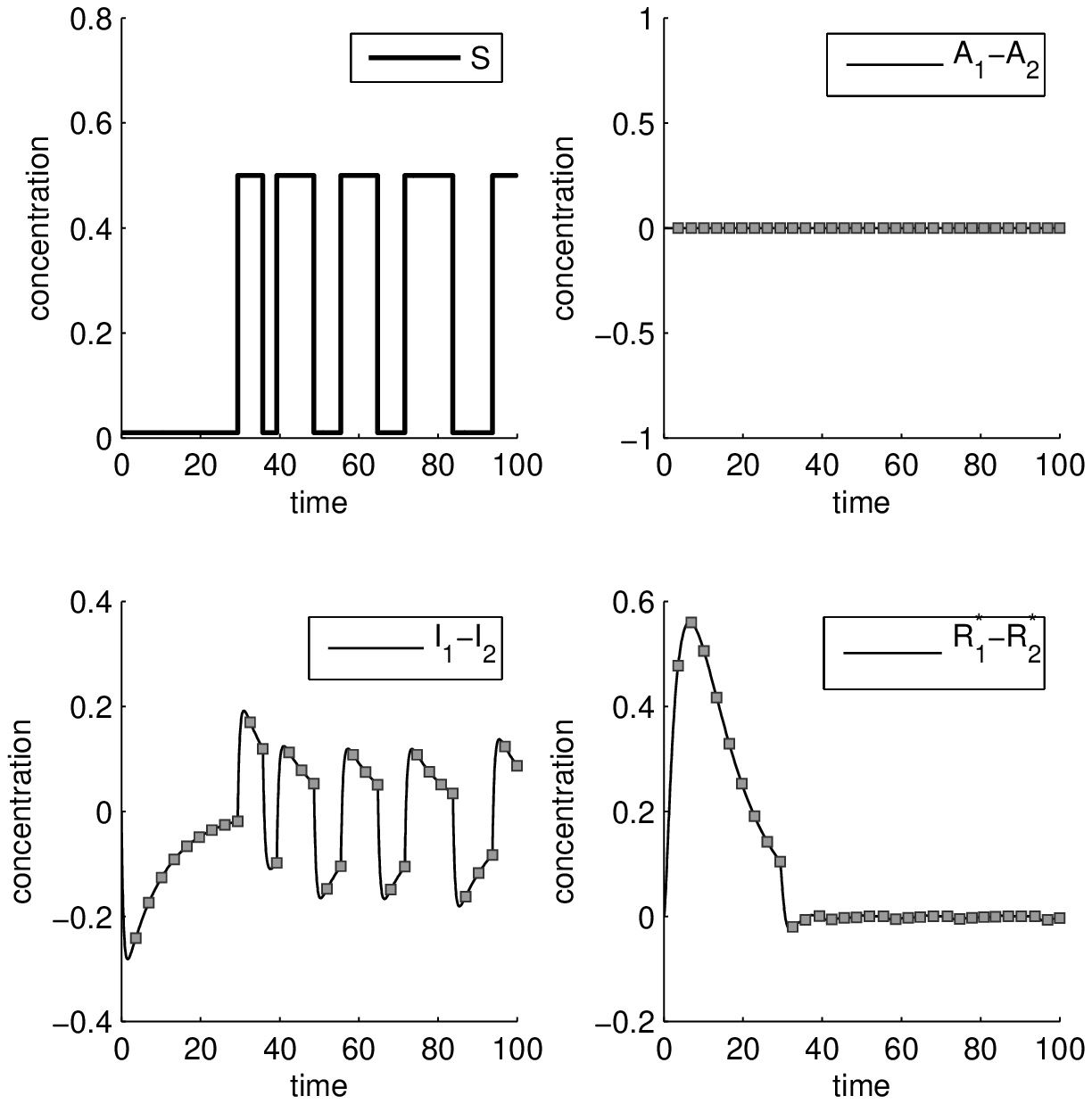}}
\subfigure{\includegraphics[width=0.49\columnwidth]{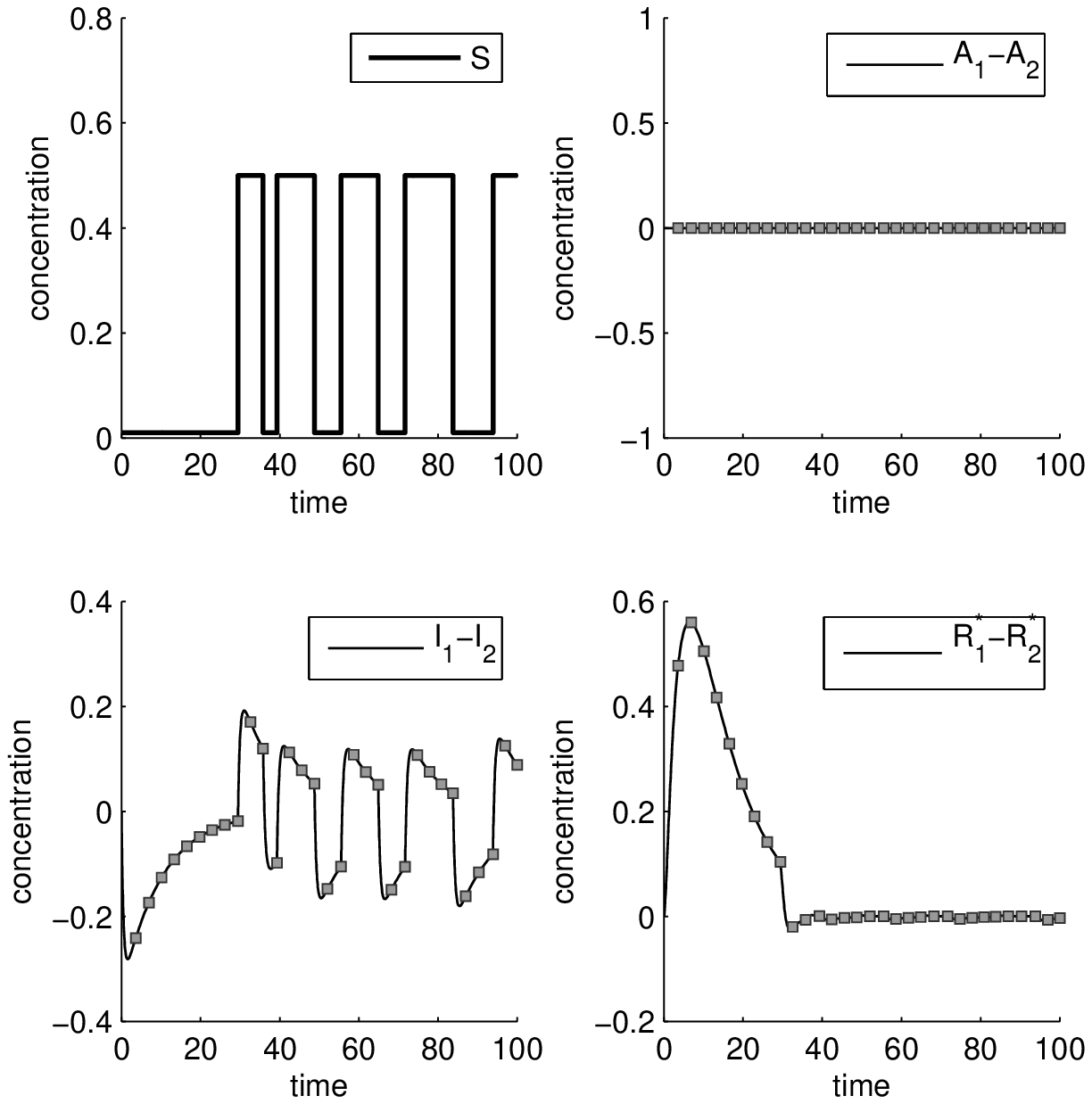}}
\end{center}
\caption{The model \corrFIVEDOM{variable} distance functions $(S_1-S_2)$, $(A_1-A_2)$, $(I_1-I_2)$ and $(R_1-R_2)$ are shown for the same solution design as for problem $\mathbf{P}_{\widetilde{\Theta}^1}$ \corrFIVEDOM{after the next robustification step} and for the final design (right) for two models describing signal sensing in dictyostelium discoideum. One square represents one measurement time point.}
\label{amoebaDESIGN2}
\end{figure}

\begin{figure}[ht]
\psfrag{label1}{\parbox{4cm}{\mbox{} \\  $N$}}
\psfrag{label2}{\parbox{4cm}{\vspace{-0.1cm}\hspace{-0.0cm} $\Delta_{RG}$}}
\begin{center}
\subfigure[with homotopy strategy]{\includegraphics[width=0.49\columnwidth]{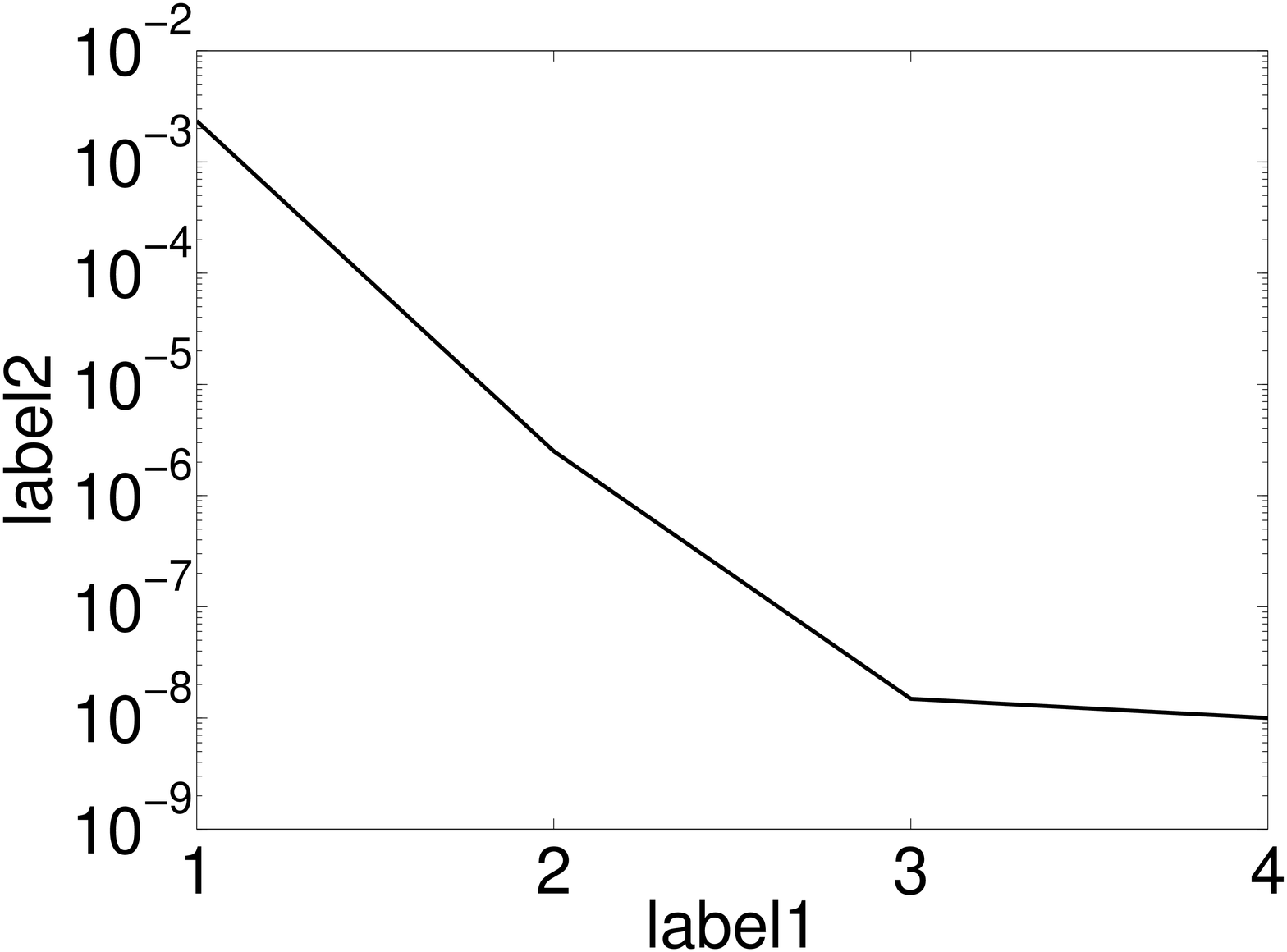}}
\subfigure[without homotopy strategy]{\includegraphics[width=0.49\columnwidth]
{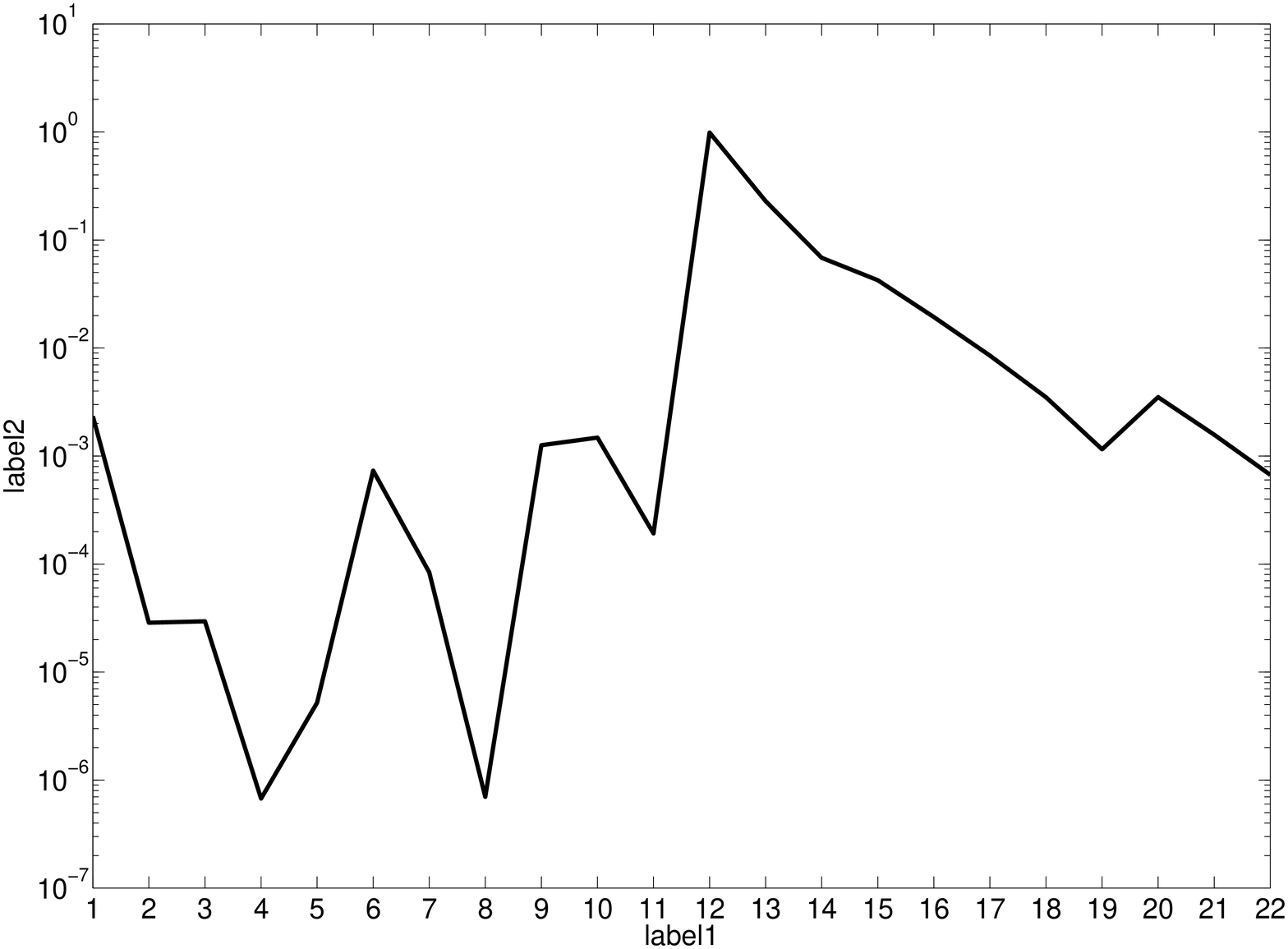}}
\end{center}
\caption{\corrFIVEDOM{T}he robustification gap $\Delta_{RG}$ is plotted versus the number of iterations $N$ of Algorithm \ref{relaxationAlgorithm} for the setup with two models describing signal sensing in dictyostelium discoideum, with homotopy strategy \corrFIVEDOM{(left)} and without homotopy strategy \corrFIVEDOM{(right)}.}
\label{amoebaDESIGN3}
\end{figure}

\begin{figure}[ht]
\begin{center}
\psfrag{label1}{\parbox{4cm}{\mbox{} \\ $N$}}
\psfrag{label2}{\parbox{4cm}{\vspace{-0.1cm}\hspace{-1.3cm} objective value of $\mathbf{P}_{\widetilde{\Theta}_N}$}}
\subfigure[with homotopy strategy]{\includegraphics[width=0.49\columnwidth]{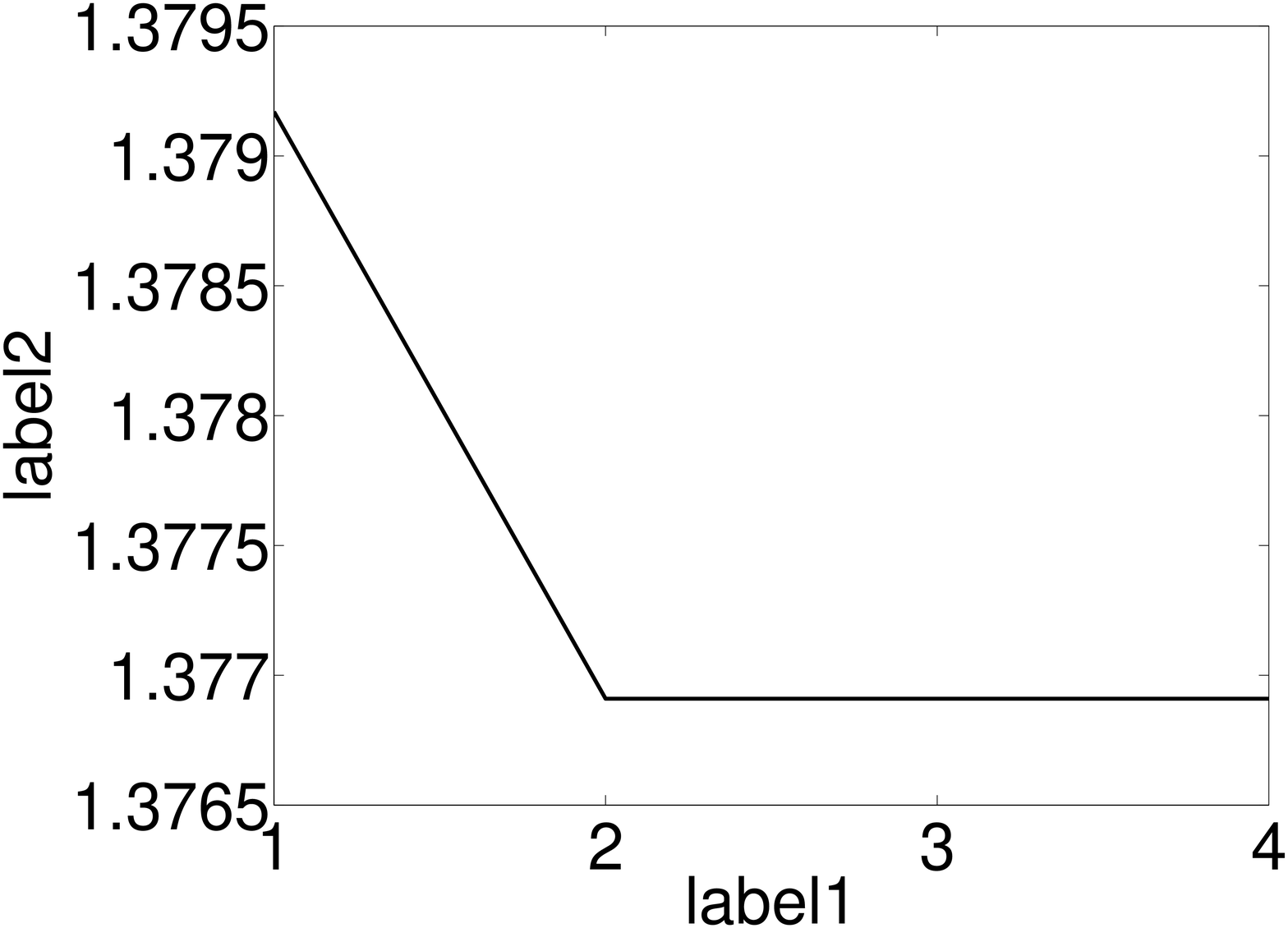}}
\subfigure[without homotopy strategy]{\includegraphics[width=0.49\columnwidth]{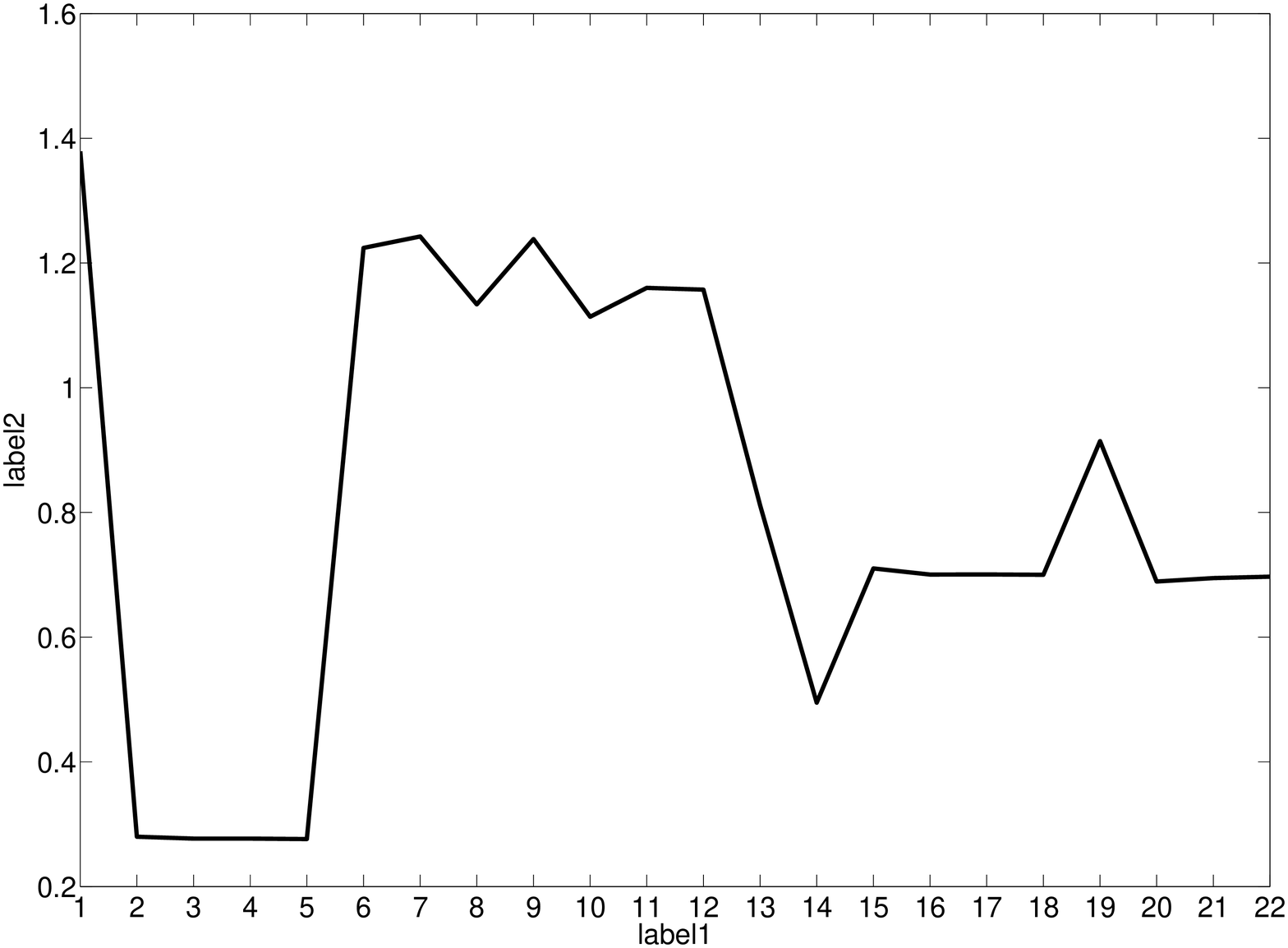}}
\end{center}
\caption{\corrFIVEDOM{T}he objective value of problem $\mathbf{P}_{\widetilde{\Theta}^N}$ is plotted versus the number of iterations $N$ of Algorithm \ref{relaxationAlgorithm} for the setup with two models describing signal sensing in dictyostelium discoideum, with homotopy strategy \corrFIVEDOM{(left)} and without homotopy strategy \corrFIVEDOM{(right)}.}
\label{amoebaDESIGN4}
\end{figure}

\begin{figure}[ht]
\begin{center}
\subfigure[with homotopy strategy]{\includegraphics[width=0.49\columnwidth]{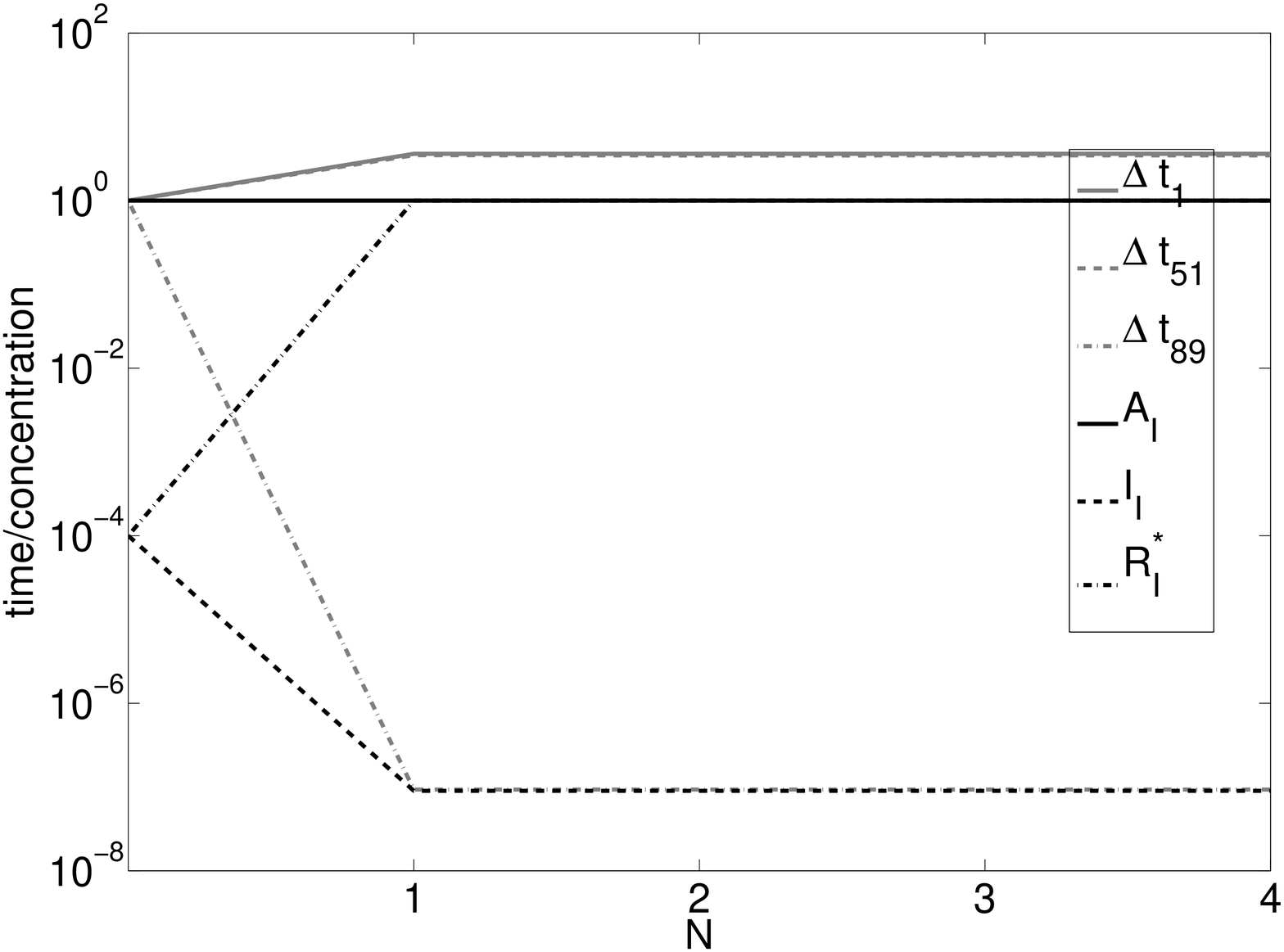}}
\subfigure[without homotopy strategy]{\includegraphics[width=0.49\columnwidth]{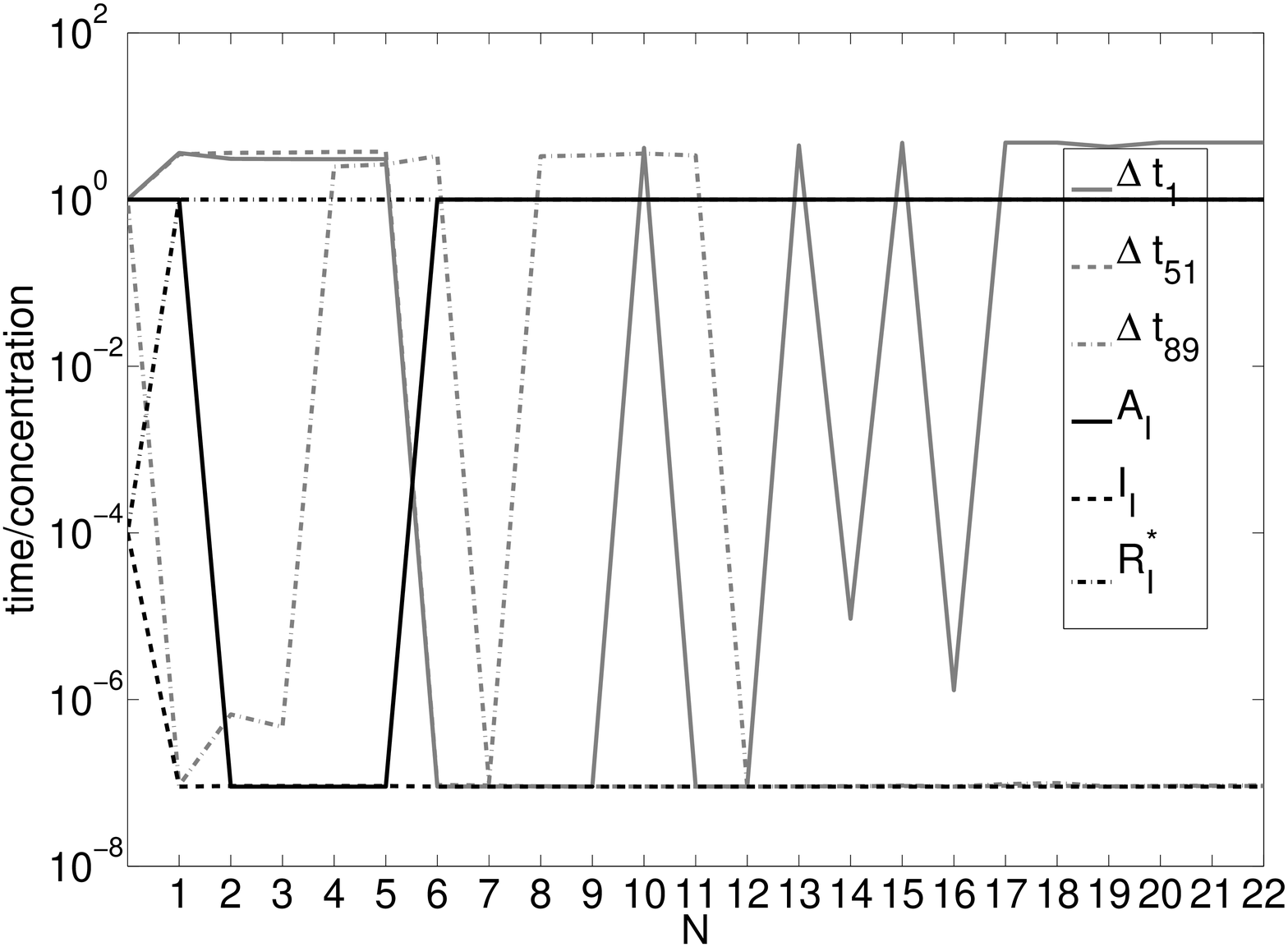}}
\end{center}
\caption{\corrFOURDOM{A selection of design variables calculated as solution of problem $\mathbf{P}_{\widetilde{\Theta}^{N}}$  is plotted versus the number of iterations $N$ of Algorithm \ref{relaxationAlgorithm} for the setup with two models describing signal sensing in dictyostelium discoideum, on the left with homotopy strategy and on the right without homotopy strategy.}}
\label{amoebaDESIGN5}
\end{figure}

The second test case is \corrFIVEDOM{the discrimination of} two models describing the chemotactic response in the amoeba \textit{Dictyostelium discoideum } as presented in \cite{Melykuti2010} using the framework presented in Section \ref{theory} and \ref{maxminTheory}. The two models describe the adaption mechanism observed when amoebae encounter the chemoattractant cAMP \cite{Levchenko2002}, see figure \ref{amoebaModel}.
For both models, a chemotaxis response regulator $R$ gets activated ($R^*$) by an activator enzyme $A$, when a cAMP ligand $S$ appears. But the deactivating mechanism determined by the interaction with an inhibitor molecule $I$ differs for both models. Both models comprise mass action kinetics in form of ODE.\\\\
In model 1 the activator enzyme as well as the inhibitor enzyme are regulated by the external signal , which is proportional to the cAMP concentration $S$. The overall model in this case is given by,
\begin{equation}
 \begin{split}
  \dot{A}_1&=-k_{-a}A_1+k_aS_1\\
  \dot{I}_1&=-k_{-i}I_1+k_{i_1}S_1\\
  \dot{R}_1^*&=-(k_rA_1+k_{-r}I_1)R_1^*+k_rR_TA_1,
 \end{split}
\end{equation}
where $k_{-a}$, $k_a$, $k_{-i}$, $k_{i_1}$, $k_r$  and $k_{-r}$ are the mass action rate constants and  $R_T:= R^*+R$ is the total amount of the response regulator.\\\\
In model 2 the inhibitory molecule $I$ is activated through the indirect action of activator $A$ instead of direct activation by sensing ligand binding. The overall model in this case is given by,
\begin{equation}
 \begin{split}
  \dot{A}_2&=-k_{-a}A_2+k_aS_2\\
  \dot{I}_2&=-k_{-i}I_2+k_{i_2}A_2\\
  \dot{R}_2^*&=-(k_rA_2+k_{-r}I_2)R_2^*+k_rR_TA_2,
 \end{split}
\end{equation}
where $k_{-a}$, $k_a$, $k_{-i}$, $k_{i_2}$, $k_r$  and $k_{-r}$ are the mass action rate constants and  $R_T:= R^*+R$ is the total amount of the response regulator.\\
For modeling details we refer to \cite{Melykuti2010}. We have extended these systems of ordinary differential equations by an additional state corresponding to the cAMP ligand $S$ with $\dot{S}=0$.
By allowing species concentration perturbations $c$ only to the state $S$ we can mimic a piecewise constant control of the system by the cAMP ligand $S$.\\\\
The experimental design parameters are the initial species concentrations of the four states namely, $A_I$, $I_I$, $R_I$, $S_I$, the measurement time points $t$ and the species concentration perturbation $c$ with respect to S. We discard the condition that either a measurement or a perturbation can be performed since in that setting by use of the perturbations $c$ we mimic a piecewise constant \corrFIVEDOM{input} control $S$ and therefore that restriction seems unnatural. Again for simplicity we consider the homoscedastic case with equal variances i.e. $v_1 = v_2 =\sigma^2$, where $\mathcal{I}(P_1:P_2,\mathcal{O}_1)$ reduces now to
\begin{equation}
 \mathcal{I}(P_1:P_2,\mathcal{O}_1)=\sum_{i=1}^n  \corrFIVEDOM{\mathcal{H}'} (t_i)\left(  (A^i_1-A^i_2)^2+(I^i_1-I^i_2)^2+(R^{*,i}_1-R^{*,i}_2)^2      \right) .
\end{equation}
The parameters $k_{-a}$, $k_a$, $k_{-i}$, $k_{i_1}$, $k_r$, $k_{-r}$ and $R_T$ are regarded as known and fixed, their values are given in Table \ref{paraValueNumExampl2}.
\begin{table}[h]
\begin{center}
\begin{tabular}{||c|c|c|c|c|c|c||}
\hline 
$k_{-a}$ & $k_a$ & $k_{-i}$ & $k_{i_1}$ & $k_r$ & $k_{-r}$ & $R_T$\tabularnewline
2.0 & 3.0 & 0.1 & 1.0 & 1.0 & 1.0 & 23/30 \tabularnewline
\hline
\end{tabular}
\end{center}
\caption{Parameter values for the fix values within model 1 and model 2.}
\label{paraValueNumExampl2}
\end{table} 
Parameter $k_{i_2}$ is regarded as unknown and subject to robustification. The range of the parameter $k_{i_2}$ is set to $k_{i_2}\in[0,2]$.\\
The optimal design is calculated within a fixed time window with $T_{\mathrm{end}}=100$. $100$ equally spaced possible measurement points are defined in the initial state of the optimization procedure. The distance vector $\Delta t$ between time points is subject to design and each entry is restricted to $\Delta t_i\in[10^{-7},10^{19}],$ $i\in \{1,...,100\}$.\\
The free perturbation vectors $c_i$, $i\in\{11,21,31,41,51,61,71,81,91\}$ are not restricted. The initial values are set to $c_i=0$, $i\in\{11,21\}$, $c_{31}=0.3$, $c_i=-0.48$, $i\in\{41,61,81\}$ and $c_i=0.48$, $i\in\{51,71,91\}$.\\
The initial species concentrations which are also subject to the experimental design are restricted 
to $S_I \in [0.01,0.5]$, $A_I \in [10^{-7},1]$, $I_I \in [10^{-7},1]$ and $R_I \in [10^{-7},I]$. The initial values are set to $S_I=0.2$, $A_I=1.0$, $I=10^{-4}$ and $R=10^{-4}$. The multiple shooting intermediate variables for the species $S$ are restricted to $s_i\in[0.01,0.5]$ to restrict the piecewise constant control to this interval. The parameters of the switching functions $\corrFIVEDOM{\mathcal{H}'}(t_i)$ are chosen as $a_1=5.0$ and $b_1=2.5$.
The algorithmic settings are summarized in Table \ref{numericSettingsExample2}.
\begin{table}[h]
\begin{center}
\begin{tabular}{||c|c|c||c|c||}
\hline 
\multicolumn{3}{||c||}{Optimization settings}  & \multicolumn{2}{c||}{Integrator settings} \tabularnewline
\hline
 $P$ & $\delta$ &  IPOPT-tol: Step 1./Step 2. & relTol/absTol & relTolSens/absTolSens \tabularnewline
 $5$ & $10^{-8}$ & $10^{-10}$/$10^{-11}$ & $10^{-14}$/$10^{-14}$ &  $10^{-14}$/$10^{-14}$ \tabularnewline
\hline
\end{tabular}
\end{center}
\caption{On the left \corrFIVEDOM{hand side} the optimization settings are listed comprising the IPOPT stopping tolerances for Step 1. and Step 2. of Algorithm \ref{relaxationAlgorithm} and on the right \corrFIVEDOM{hand side} the integration tolerances for the nominal trajectory and the first \corrFIVEDOM{order} sensitivities. We use the IPOPT option  ``honor\_original\_bounds=no'' for Step 1. and Step 2. of Algorithm \ref{relaxationAlgorithm}.}
\label{numericSettingsExample2}
\end{table}\\
With these design conditions we start the optimization procedure twice. First by use of the homotopy strategy for successive problems $\mathbf{P}_{\widetilde{\Theta}^{N+1}}$ with $10$ homotopy steps.\\
Since the ``discriminating power'' of the experimental setup is very low \corrFIVEDOM{in this case}, i.e. the deviation between the two models is small, we plot the distance functions $(S_1-S_2)$, $(A_1-A_2)$, $(I_1-I_2)$ and $(R_1-R_2)$ for the initial state and for the solution of problem $\mathbf{P}_{\widetilde{\Theta}^1}$ in Figure \ref{amoebaDESIGN1}. A plot for the same functions with the same solution design as for problem $\mathbf{P}_{\widetilde{\Theta}^1}$ \corrFIVEDOM{after the next robustification step} is shown in Figure \ref{amoebaDESIGN2}. The final design is also shown in Figure \ref{amoebaDESIGN2}.\\
A plot of the robustification gap $\Delta_{RG}$ for each iteration $N$ of Algorithm \ref{relaxationAlgorithm} is shown in Figure \ref{amoebaDESIGN3} (left). A plot of the objective value of problem $\mathbf{P}_{\widetilde{\Theta}^{N}}$ for each iteration $N$ of Algorithm \ref{relaxationAlgorithm} is shown in Figure \ref{amoebaDESIGN4} (left). A selection of design variables as solution\corrFIVEDOM{s} of problem $\mathbf{P}_{\widetilde{\Theta}^{N}}$ is shown in Figure \ref{amoebaDESIGN5} (left).\\\\
Second we calculate the design without the homotopy strategy. We experience huge jumps in the final objective value of problem $\mathbf{P}_{\widetilde{\Theta}^N}$ for subsequent iterations $N$ of Algorithm \ref{relaxationAlgorithm}. This is due to the fact that the final design of the former problem $\mathbf{P}_{\widetilde{\Theta}^N}$ is an infeasible starting point for the successive problem $\mathbf{P}_{\widetilde{\Theta}^{N+1}}$ in the interior point solution strategy. First the optimizer tries to force the iterates back into the feasible region and afterwards the new central path leads to a different design.\\
For this case a plot of the robustification gap $\Delta_{RG}$ for each iteration $N$ of Algorithm \ref{relaxationAlgorithm} is shown in Figure \ref{amoebaDESIGN3} (right). A plot of the objective value of problem $\mathbf{P}_{\widetilde{\Theta}^{N}}$ for each iteration $N$ of Algorithm \ref{relaxationAlgorithm} is shown in Figure \ref{amoebaDESIGN4} (right). A selection of design variables as solution\corrFIVEDOM{s} of problem $\mathbf{P}_{\widetilde{\Theta}^{N}}$ is shown in Figure \ref{amoebaDESIGN5} (right).
As one can clearly see\corrFIVEDOM{,} the homotopy strategy helps to considerably stabilize Algorithm \ref{relaxationAlgorithm}.
\section{Conclusion}
We present a framework for the robust computation of optimal experimental
designs for the purpose of model discrimination. The theoretical framework as
well as the numerical realization by utilization of an outer approximation
algorithm are worked out. A strategy for the numerical stabilization of the
algorithm by use of a homotopy approach is suggested. The optimization
procedure is successfully exemplified on two biological model systems. In our
examples we clearly found that the homotopy approach is significantly superior
to \corrFIVEDOM{a cold} start of successive design problems
$\mathbf{P}_{\widetilde{\Theta}^{N+1}}$. For the first test case, the
discrimination of two models describing glycolytic oscillations, the outer
approximation scheme completely fails to reach the desired accuracy $\delta$
without homotopy strategy. For the second test case, the discrimination of two
models describing signal sensing in dictyostelium discoideum, the outer
approximation scheme also fails without warmstart, however the homotopy
strategy also works with only two homotopy steps (not presented in this
paper). This indicates the need of a step size strategy for reasons of
efficiency which will be a next step in our work.
\section{Acknowledgement}
The authors thank the anonymous reviewers for helpful comments and suggestions.\\\\
The authors gratefully acknowledge the Freiburg Initiative for Systems Biology
(FRISYS), part of the BMBF FORSYS systems biology initiative, the 
Freiburg excellence cluster Centre for Biological Signalling Studies
(BIOSS), the Helmholtz alliance Systems Biology of Cancer and the Nephage
iniative (BMBF Gerontosys II) for various support and funding.

\bibliography{/home/dskanda/literature/literature}

\begin{thebibliography}{10}

\bibitem{Albersmeyer2010}
Jan Albersmeyer.
\newblock {\em Adjoint-based algorithms and numerical methods for sensitivity
  generation and optimization of large scale dynamic systems}.
\newblock PhD thesis, University of Heidelberg, Heidelberg, December 2010.

\bibitem{Albersmeyer2008a}
Jan Albersmeyer and Hans~Georg Bock.
\newblock Sensitivity generation in an adaptive {BDF}-method.
\newblock In {\em Modeling, Simulation and Optimization of Complex Processes:
  Proceedings of the Third International Conference on High Performance
  Scientific Computing}. Springer, 2008.

\bibitem{Apgar2008}
Joshua~F Apgar, Jared~E Toettcher, Drew Endy, Forest~M White, and Bruce Tidor.
\newblock Stimulus design for model selection and validation in cell signaling.
\newblock {\em PLoS Computational Biology}, 4(2):e30, 02 2008.

\bibitem{Arora1995}
J.~S. Arora, O.~A. Elwakeil, A.~I. Chahande, and C.~C. Hsieh.
\newblock Global optimization methods for engineering applications: A review.
\newblock {\em Structural and Multidisciplinary Optimization}, 9:137--159,
  1995.
\newblock 10.1007/BF01743964.

\bibitem{ATKINSON1975}
A.~C. Atkinson and V.~V. Fedorov.
\newblock The design of experiments for discriminating between two rival
  models.
\newblock {\em Biometrika}, 62(1):57--70, 1975.

\bibitem{Balsao2008}
E.~Balsa-Canto, A.~A. Alonso, and J.~R. Banga.
\newblock Computational procedures for optimal experimental design in
  biological systems.
\newblock {\em IET Systems Biology}, 2(4):163--172, July 2008.

\bibitem{Bauer2000}
I.~Bauer, H.~G. Bock, S.~K{\"o}rkel, and J.~P. Schl{\"o}der.
\newblock Numerical methods for optimum experimental design in {DAE} systems.
\newblock {\em Journal of Computational and Applied Mathematics}, 120:1--25,
  2000.

\bibitem{Bell2010}
Bradley~M. Bell.
\newblock Automatic differentiation software cppad., 2010.

\bibitem{Bell2008}
Bradley~M. Bell and James~V. Burke.
\newblock Algorithmic differentiation of implicit functions and optimal values.
\newblock In Christian~H. Bischof, H.~Martin B{\"u}cker, Paul~D. Hovland, Uwe
  Naumann, and J.~Utke, editors, {\em Advances in Automatic Differentiation},
  pages 67--77. Springer, Berlin, 2008.

\bibitem{Bernacki2009}
Joseph~P. Bernacki and Regina~M. Murphy.
\newblock Model discrimination and mechanistic interpretation of kinetic data
  in protein aggregation studies.
\newblock {\em Biophysical Journal}, 96:2871--2887, 2009.

\bibitem{Biegler2001}
Lorenz~T. Biegler, Arturo~M. Cervantes, and Andreas W{\"a}chter.
\newblock Advances in simultaneous strategies for dynamic process optimization.
\newblock {\em Optimization, Chemical Engineering Science}, 57:575--593, 2001.

\bibitem{Billingsley1986}
Patrick Billingsley.
\newblock {\em Probability and Measure}.
\newblock John Wiley \& Sons Inc, 1986.

\bibitem{Bock1987}
Hans~Georg Bock.
\newblock {R}andwertproblemmethoden zur {P}arameteridentifizierung in
  {S}ystemen nichtlinearer {D}ifferentialgleichungen.
\newblock In {\em Bonner Mathematische Schriften}, volume 183. University of
  Bonn, 1987.

\bibitem{Bock1984}
Hans~Georg Bock and Karl~J. Plitt.
\newblock A multiple shooting algorithm for direct solution of optimal control
  problems.
\newblock In {\em Proceedings of the Ninth IFAC World Congress, Budapest}.
  Pergamon, Oxford, 1984.

\bibitem{Burnham2002}
Kenneth~P. Burnham and David~R. Anderson.
\newblock {\em Model Selection and Multimodel inference: A practical
  information-theoretic approach}.
\newblock Springer, 2002.

\bibitem{Byrne1975}
G.~D. Byrne and A.~C. Hindmarsh.
\newblock A polyalgorithm for the numerical solution of ordinary differential
  equations.
\newblock {\em {ACM} Transactions on Mathematical Software}, 1(1):71--96, 1975.

\bibitem{Calvo1993}
M.~Calvo, J.~I. Montijano, and L.~R{\'a}ndez.
\newblock On the change of step size in multistep codes.
\newblock {\em Numerical Algorithms}, 4:283--304, 1993.

\bibitem{Chernoff1956}
Herman Chernoff.
\newblock Large-sample theory: Parametric case.
\newblock {\em The Annals of Mathematical Statistics}, 27(1):pp. 1--22, 1956.

\bibitem{Cooney1995}
M.~J. Cooney and K.~A. McDonald.
\newblock Optimal dynamic experiments for bioreactor model discrimination.
\newblock {\em Applied Microbiology and Biotechnology}, 43:826--837, 1995.

\bibitem{Goldbeter1996}
Albert Goldbeter.
\newblock {\em Biochemical oscillations and cellular rhythms: The molecular
  bases of periodic and chaotic behaviour}.
\newblock Cambridge University Press, 1996.

\bibitem{Hettich1993}
R.~Hettich and K.~O. Kortanek.
\newblock Semi-infinite programming: Theory, methods, and applications.
\newblock {\em SIAM Review}, 35(3):pp. 380--429, 1993.

\bibitem{Horn1987}
R.~Horn.
\newblock Statistical methods for model discrimination. applications to gating
  kinetics and permeation of the acetylcholine receptor channel.
\newblock {\em Biophysical Journal}, 51:255--263, 1987.

\bibitem{HSL}
HSL.
\newblock A collection of fortran codes for large-scale scientific computation.
  {S}ee http://www.hsl.rl.ac.uk, 2007.

\bibitem{RishiJain2006}
Rishi Jain, Andrea~L. Knorr, Joseph Bernacki, and Ranjan Srivastava.
\newblock Investigation of bacteriophage ms2 viral dynamics using model
  discrimination analysis and the implications for phage therapy.
\newblock {\em Biotechnology Progress}, 22(6):1650--1658, 2006.

\bibitem{Korkel1999}
S.~K{\"o}rkel, I.~Bauer, H.~G. Bock, and J.~P. Schl{\"o}der.
\newblock A sequential approach for nonlinear optimum experimental design in
  {DAE} systems.
\newblock In F.~Keil, W.~Mackens, H.~Voss, , and J.~Werther, editors, {\em
  Scientific Computing in Chemical Engineering II}, volume~2. Springer Verlag,
  Berlin, 1999.

\bibitem{Kremling2004}
A.~Kremling, S.~Fischer, K.~Gadkar, F.~J. Doyle, T.~Sauter, E.~Bullinger,
  F.~Allg\"{o}wer, and E.~D. Gilles.
\newblock A benchmark for methods in reverse engineering and model
  discrimination: problem formulation and solutions.
\newblock {\em Genome Research}, 14(9):1773--1785, September 2004.

\bibitem{Kreutz2009}
Clemens Kreutz and Jens Timmer.
\newblock Systems biology: experimental design.
\newblock {\em FEBS Journal}, 276:923--942, 2009.

\bibitem{Kullback1997}
Solomon Kullback.
\newblock {\em Information Theory and Statistics}.
\newblock Dover Publications Inc., 1997.

\bibitem{Lacey1984}
Laurence Lacey and Adrian Dunne.
\newblock The design of pharmacokinetic experiments for model discrimination.
\newblock {\em Journal of Pharmacokinetics and Pharmacodynamics}, 12:351--365,
  1984.

\bibitem{Leineweber1998}
Daniel~B. Leineweber.
\newblock {\em Efficient Reduced SQP Methods for the Optimization of Chemical
  Processes Described by Large Sparse DAE Models}.
\newblock PhD thesis, University of Heidelberg, 1998.

\bibitem{Levchenko2002}
A~Levchenko and PA~Iglesias.
\newblock Models of eukaryotic gradient sensing: Application to chemotaxis of
  amoebae and neutrophils.
\newblock {\em Biophysical Journal}, 82:50--63, 2002.

\bibitem{Lopez-Fidalgo2007}
J.~L{\'o}pez-Fidalgo, C.~Tommasi, and P.~C. Trandafir.
\newblock An optimal experimental design criterion for discriminating between
  non-normal models.
\newblock {\em Journal of the Royal Statistical Society Series B},
  69(2):231--242, 2007.

\bibitem{Melykuti2010}
Bence Melykuti, Elias August, Antonis Papachristodoulou, and Hana El-Samad.
\newblock Discriminating between rival biochemical network models: three
  approaches to optimal experiment design.
\newblock {\em BMC Systems Biology}, 4(1):38, 2010.

\bibitem{Myung2009}
Jay~I. Myung and Mark~A. Pitt.
\newblock Optimal experimental design for model discrimination.
\newblock {\em Psychological review}, 116(3):499--518, July 2009.

\bibitem{Polak1969}
E.~Polak.
\newblock On the convergence of optimization algorithms.
\newblock {\em Rev. Fran\c caise Informat. Recherche Op\'erationnelle},
  3(16):17--34, 1969.

\bibitem{Polak1987}
E.~Polak.
\newblock On the mathematical foundations of nondifferentiable optimization in
  engineering design.
\newblock {\em SIAM Review}, 29(1):pp. 21--89, 1987.

\bibitem{Polak1993}
E.~Polak.
\newblock On the use of consistent approximations in the solution of
  semi-infinite optimization and optimal control problems.
\newblock {\em Mathematical Programming}, 62:385--414, 1993.
\newblock 10.1007/BF01585175.

\bibitem{Polak1997}
Elijah Polak.
\newblock {\em Optimization: Algorithms and Consistent Approximations}.
\newblock Springer, 1997.

\bibitem{Pronzato1988}
Luc Pronzato and Eric Walter.
\newblock Robust experiment design via maximin optimization.
\newblock {\em Mathematical Biosciences}, 89(2):161 -- 176, 1988.

\bibitem{Perez2009}
Victor Pérez, John Renaud, and Layne Watson.
\newblock Homotopy curve tracking in approximate interior point optimization.
\newblock {\em Optimization and Engineering}, 10:91--108, 2009.
\newblock 10.1007/s11081-008-9042-6.

\bibitem{Salmon1968}
D.~Salmon.
\newblock Minimax controller design.
\newblock {\em Automatic Control, IEEE Transactions on}, 13(4):369 -- 376, aug.
  1968.

\bibitem{Shimizu1980}
Kiyotaka Shimizu and Eitaro Aiyoshi.
\newblock Necessary conditions for min-max problems and algorithms by a
  relaxation procedure.
\newblock {\em IEEE Transactions on Automatic Control}, 25(1):62--66, 1980.

\bibitem{Skanda2010}
Dominik Skanda and Dirk Lebiedz.
\newblock An optimal experimental design approach to model discrimination in
  dynamic biochemical systems.
\newblock {\em Bioinformatics}, 26(7):939--945, 2010.

\bibitem{Stewart1998}
W.~E. Stewart, Y.~Shon, and G.~E.~P. Box.
\newblock Discrimination and goodness of fit of multiresponse mechanistic
  models.
\newblock {\em AIChE Journal}, 44(6):1404--1412, 1998.

\bibitem{Stoer2002}
Josef Stoer and Roland Bulirsch.
\newblock {\em Introduction to Numerical Analysis}.
\newblock Number~12 in Texts in Applied Mathematics. Springer, New~York, third
  edition, 2002.

\bibitem{Stricker1994}
C.~Stricker, S.~Redman, and D.~Daley.
\newblock Statistical analysis of synaptic transmission: model discrimination
  and confidence limits.
\newblock {\em Biophysical Journal Of The Royal Statistical Society Series B},
  67:532--547, 1994.

\bibitem{Takors1997}
R.~Takors, W.~Wiechert, and D.~Weuster-Botz.
\newblock Experimental design for the identification of macrokinetic models and
  model discrimination.
\newblock {\em Biotechnol Bioeng}, 56(5):564--576, Dec 1997.

\bibitem{Timmer2004}
Jens Timmer, T.~G. M{\"u}ller, I.~Swameye, O.~Sandra, and U.~Klingm{\"u}ller.
\newblock Modeling the nonlinear dynamics of cellular signal transduction.
\newblock {\em International Journal of Bifurcation and Chaos},
  14(6):2069--2079, 2004.

\bibitem{Ucinski2004}
D.~Uci{\'n}ski and B.~Bogacka.
\newblock {T}-optimum designs for multiresponse dynamic heteroscedastic models.
\newblock In A.~Di Bucchianico and H.~Lauter, editors, {\em Proc. of the 7th
  International Workshop on Model-Oriented Design and Analysis}, pages
  191--199. Physica Verlag, 2004.

\bibitem{Waechter2002}
Andreas W{\"a}chter.
\newblock {\em An Interior Point Algorithm for Large-Scale Nonlinear
  Optimization with Applications in Process Engineering}.
\newblock PhD thesis, Carnegie Mellon University, 2002.

\bibitem{Waechter2006}
Andreas W{\"a}chter and Lorenz~T. Biegler.
\newblock On the implementation of a primal-dual interior point filter line
  search algorithm for large-scale nonlinear programming.
\newblock {\em Mathematical Programming}, 106(1):25--57, 2006.

\end{thebibliography}
\bibliographystyle{plain}
\end{document}